\documentclass[pdftex,a4paper,leqno,numbers=noenddot]{scrartcl}
\usepackage[latin1]{inputenc}
\usepackage[T1]{fontenc}
\usepackage{amssymb}
\usepackage{mathrsfs}
\usepackage{amsmath}
\usepackage{graphicx}
\usepackage{array}
\usepackage{fancyhdr}
\usepackage[format=hang,
						font={footnotesize,sf},
						labelfont={bf},
						margin=1cm,
						aboveskip=5pt,
						position=bottom]{caption}
\usepackage[numbers]{natbib}
\numberwithin {equation} {section}
\numberwithin {table} {section}
\begin{document}
\thispagestyle{empty}
\begin{center}
\textbf{\Large{The Asymptotic Covariance Matrix of the Odds Ratio Parameter Estimator in Semiparametric Log-bilinear Odds Ratio Models}}\par\bigskip
\textbf{\large{Angelika Franke, Gerhard Osius}}\par\bigskip
\large{Faculty 3 / Mathematics / Computer Science}\\\large{University of Bremen}
\end{center}
\hrule

\def\a{\alpha}
\def\b{\beta}
\def\r{\rho}
\def\g{\gamma}
\def\l{\lambda}
\def\D{\Delta}
\def\ve{\VaRepsilon}
\def\s{\sigma}
\def\k{\kappa}
\def\t{\theta}
\def\VaR{\VaRsigma}
\def\O{\Omega}
\def\RX{\mathbb{R}^{M_X}}
\def\RY{\mathbb{R}^{M_Y}}
\def\RhX{\mathbb{R}^{L_X}}
\def\RhY{\mathbb{R}^{L_Y}}
\newtheorem{Th}{Theorem}
\newtheorem{Co}{Corollary}

\section*{Abstract}
The association between two random variables is often of primary interest in statistical research. In this paper semiparametric models for the association between random vectors X and Y are considered which leave the marginal distributions arbitrary. Given that the odds ratio function comprises the whole information about the association the focus is on bilinear log-odds ratio models and in particular on the odds ratio parameter vector $\t$. The covariance structure of the maximum likelihood estimator $\hat{\t}$ of $\t$ is of major importance
for asymptotic inference. To this end different representations of the estimated covariance matrix are derived for conditional and unconditional sampling schemes and different asymptotic approaches depending on whether X and/or Y has finite or arbitrary support. The main result is the invariance of the estimated asymptotic covariance matrix of $\hat{\t}$ with respect to all above approaches. As applications we compute the asymptotic power for tests of linear hypotheses about $\t$---with emphasis to logistic and linear regression models---which allows to determine the necessary sample size to achieve a wanted power.\par\bigskip\noindent
\small{\textbf{AMS 2000 subject classifications:} 
Primary 62F12; secondary 62D05, 62H17, 62J05, 62J12.\par\bigskip\noindent
\textbf{Keywords:}
Odds ratio; asymptotic; covariance matrix; conditional sampling; semiparametric; log-linear models; log-bilinear association; logistic regression; linear regression.}
\par\bigskip\noindent
\hrule
\pagenumbering{arabic}

\fancyhf{}
\fancyhead[L]{\nouppercase{\leftmark}}
\fancyfoot[C]{\thepage}
\setcounter{page}{1}
\section{Introduction and Outline}
The question how a random output vector $Y$ of a system (e.g. the health status of a human) is associated to a random input vector $X$ (e.g. consumption of tobacco and alcohol, environmental pollution and other risk factors) is of major importance in
statistical science. If the association between $X$ and $Y$ is of primary interest, then a semi-parametric association is appropriate which leaves the marginal distributions of $X$ and $Y$ arbitrary. However, the association is completely determined by the odds-ratio function $OR(x, y)$ for the joint density $p(x, y)$ with respect to fixed reference values $x_0$ and $y_0$ (cf.\ Osius 2004, 2009 \cite{Met, Anals}):
\begin{equation}
OR(x, y) =\frac{p(x, y) \cdot p(x_0, y_0)}{p(x, y_0) \cdot p(x_0, y)}.
\end{equation}
A semi-parametric odds-ratio model specifies this function up to an unknown parameter vector $\t$, but leaves marginal distributions arbitrary. An important class are log-bilinear odds-ratio models given by
\begin{equation}
\log OR(x, y) = \tilde{x}^T\t \tilde{y}
\end{equation}
where $\tilde{x}$ and $\tilde{y}$ are known vector-valued functions of $x$ and $y$ which may coincide with $x$ and $y$, respectively. The association structure of some widely used regression models is log-bilinear, e.g. generalized linear models with canonical link (for univariate $Y$), multivariate linear logistic regression (for $Y$ with finite support) and multivariate
linear regression. An advantage of odds-ratio models over these regression models is that inference about the association parameter $\t$ may also be obtained from samples drawn conditionally on $Y$ (instead of $X$). Generalizing an important result by Prentice and Pyke, 1979 \cite{PP}, it has been shown in Osius, 2009 \cite{Anals}, that the estimator $\hat{\t}$ and its estimated asymptotic covariance matrix $Cov_{\infty}(\hat{\t})$ for samples conditional on $Y$ are exactly the same as if the sample had been drawn conditionally on $X$. The purpose of this paper is to derive different representations of this covariance
matrix on which statistical analysis (e.g. tests and confidence regions) are based. These results are applied to compute the asymptotic power for tests of linear hypothesis about $\t$ which allows to determine the sample sizes necessary to achieve a wanted
power.\par\bigskip
\noindent
A given random sample $(X_i, Y_i),\ i = 1,\;\ldots,\;n$ containing $J+1$ different $X$-values $X_{(0)},\;\ldots,\;X_{(J)}$ and $K+1$ different $Y$-values $Y_{(0)},\;\ldots,\;Y_{(K)}$ can be summarized by the counts
\begin{equation}
R_{jk} = \{ i\ |\ X_i = X_{(j)},\ Y_i = Y_{(k)}\}
\end{equation}
for the observed combinations $(j,\ k)$. Although the distribution of the table $(R_{jk})$ depends on the sampling scheme (e.g. conditional on $X$ or $Y$), we will show that the estimated asymptotic covariance matrix of $\hat{\t}$ is invariant against common sampling schemes and asymptotic approaches. However we do not establish original asymptotic results here but---using mainly matrix algebra---derive different representations for asymptotic covariance matrices and in particular for $Cov_{\infty}(\hat{\t})$.\par\bigskip
\noindent
The paper is organized as follows. Section 2 gives a brief introduction to odds ratio models with emphasis on multivariate linear logistic regression (where $Y$ has finite support) and log-linear models for contingency tables (where the support of $X$ is finite too). The next section 3 deals with estimation of $\t$ under different sampling schemes (unconditional and conditional on $X$ and $Y$, respectively). Our main results are contained in section 4. Based on the work of Haberman, 1974 \cite{Hab} we first show that for contingency tables (i.e. both $X$ and $Y$ have finite support) the
asymptotic distribution of $\hat{\t}$ is invariant under the common sampling schemes and provide different representations of $Cov_{\infty}(\hat{\t})$. Looking more generally at the multivariate linear logistic regression model (with arbitrary support of $X$) and sampling conditional on $X$ we observe, that the estimated asymptotic covariance matrix $Cov_{\infty}(\hat{\t})$ is the same as for contingency tables (where $X$ has finite support). The general case allowing arbitrary supports for $X$ and $Y$ is dealt with in section 5. For sampling conditional on $Y$ and a fixed set of conditioning values we conclude that the matrix of $Cov_{\infty}(\hat{\t})$ is the same as before where both $X$ and $Y$ had finite support. As a first application we show in section 6 how our results can be used to compute the asymptotic power for testing a linear hypothesis $Q\t = 0$ and how to determine the necessary sample size to achieve a given power for a value $\t'$ of interest under the alternative. Finally we demonstrate for univariate $Y$ how the linear resp. log-linear model emerges from an odds-ratio-model by imposing additional assumptions on the conditional distribution of $Y$ (given $X$) and conclude with a short discussion of our results. The appendix contains the proofs and some results from linear algebra.

\section{Odds Ratio Models}\label{sec2}
Consider a pair $(X,Y)$ of random vectors defined on some probability space taking values in $\O=\O_X\times\O_Y\subset \RX\times\RY$ with joint distribution $P$ and marginal distributions $P^X$ and $P^Y$. To avoid trivialities we assume that $\O_X$ and $\O_Y$ both have more than one element. Let $\nu_X$ and $\nu_Y$ be two fixed $\sigma$-finite measures on $\RX$ and $\RY$ such that $P$ has a positive density $p$ on $\O$ with respect to the product measure $\nu=\nu_X\times\nu_Y$---typically a product of Lebesgue or counting measures. The log-density can be parametrized as
\begin{equation}\label{Y5}
\log p(x,y)=\ \alpha+\rho(x)+\g(y)+\psi_{\t}(x,y),\hspace{1cm} x\in\O_X,\ y\in\O_Y
\end{equation}
with integrable functions $\r$, $\g$, $\psi$, an unknown parameter $\t\in\Theta$, and an integration constant $\a$ determined by $\int p\;d\nu=1$. 
To guarantee identifiability we assume the constraints
\begin{equation}\label{constr2.3}
\rho(x_0)=\ \g(y_0)=\ 0
\end{equation}
where $x_0\in\O_X$ and $y_0\in\O_Y$ are the reference values of the odds ratio function. 
The conditional distribution of $Y$ given $X$ has a positive density $p(y|X=x)$ given by
\begin{equation}\label{Y6}
\log p(y|X=x)=\ \g(y)+\psi_{\t}(x,y) -\delta_{\t}(x)
\end{equation}
with an integration constant $\delta_{\t}(x)$ and similarly
\begin{equation}\label{Y7}
\begin{split}
\log p(x|Y=y)&=\ \rho(x)+\psi_{\t}(x,y) -\varepsilon_{\t}(y).
\end{split}
\end{equation}
An important class of parametric association models are log-bilinear association models with respect to the transformed variables $\tilde{x} = h_X(x)$ and $\tilde{y} = h_Y (y)$ given by measurable maps $h_X:\ \RX\ \rightarrow\ \RhX$ and $h_Y:\ \RY\ \rightarrow\ \RhY$ which will always be chosen here such that $\tilde{x}_0=h_X(x_0)=0$ and $\tilde{y}_0=h_Y(y_0)=\ 0$. The functions $h_X$ and $h_Y$ are typically injective (one-to-one) but to avoid trivialities we merely assume that they are not constant. The parameter $\t$ is a $L_X\times L_Y$-matrix and the log-odds ratio function is bilinear in the transformed variables $\tilde{x}$ and $\tilde{y}$ 
\begin{equation}\label{C66}
\psi_{\t}(x,y)=\ \tilde{x}^T\t \tilde{y}
\hspace{2cm}\text{for all $x$, $y$.}
\end{equation}
This model is semiparametric in the sense that it does not restrict the marginal distributions $P^X$ and $P^Y$ except for reasonable moment conditions. More precisely, it has been shown by Osius, 2009 \cite[Sec. 3]{Anals}, that given the marginal distributions $P^X$ and $P^Y$, there exists for any $L_X\times L_Y$ matrix $\t$ a unique joint distribution $P$ with these marginals such that \eqref{C66} holds---provided the expectations $\mathbb{E}(||h_X(X)||^2)$ and $\mathbb{E}(||h_Y(Y)||^2)$ are finite, i.e the covariance matrices of $h_X(X)$ and $h_Y(Y)$ exist, and this will be assumed throughout the paper.
\par\bigskip\noindent
It will be convenient to interpret a $m\times n$ matrix $A$ as a vector $\vec{A}$ of length $mn$ obtained by placing the columns of $A$ one after another. 
Using the Kronecker product $\tilde{y}\otimes \tilde{x}$ (cf.\ appendix \ref{A.2}) the model \eqref{C66} may be rewritten as 
\begin{equation}\label{C6}
\psi_{\t}(x,y)=\ (\tilde{y}\otimes \tilde{x})^T \vec{\t} \hspace{2cm}\text{for all $x$, $y$.}
\end{equation}
Any submodel specified by a linear restriction of the form $\t = A^T \t^*B$ with given matrices $A$, $B$ and parameter matrix $\t^*$ yields a log-bilinear association too, with respect to $h^*_X = Ah_X$, $h^*_Y = Bh_Y$.\par\bigskip
\noindent
The following examples reveal that the association structure of some widely used regression models is in fact log-bilinear.
\par\bigskip
\noindent
\textbf{Example 1:} Generalized linear models\par\medskip\noindent
Let $Y$ be a univariate random variable and suppose that the conditional density of $Y$ given $X = x$ belongs to the exponential family 
\begin{equation}\label{(2.9)}
p(y | X = x) =\exp\{\phi^{-1} [y \cdot\tau(x) - b(\tau(x))] + c(y,\phi)\} 
\end{equation}
with suitable functions $b, c, \tau$ and a dispersion parameter $\phi$; compare McCullagh and Nelder (1989). Then the log-odds ratio function has
the form 
\begin{equation}
\psi(x, y) = \phi^{-1} [\tau(x)-\tau(x_0)]\cdot[y-y_0]
\end{equation}
and $\tau(x)$ is a strictly monotone function of the conditional expectation $\mu(x) = \mathbb{E}(Y | X = x)=b'(\tau(x))$. A generalized linear model with canonical link specifies the canonical parameter
\begin{equation}\label{Y2}
\tau(x) = \alpha +\tilde{x}^T \beta,
\end{equation}
where $\tilde{x}\in\RhX$ is a known vector of formal covariates and $\alpha\in\mathbb{R}$, $\beta\in\RhX$ are unknown parameters. The corresponding log-odds ratio function 
\begin{equation}\label{Y1}
\psi(x, y) = \tilde{x}^T\t y
\end{equation}
is of the form \eqref{C6} with $\tilde{y} = y$ and parameter $\t= \phi^{-1}\beta$. Note that the intercept $\alpha$ is no longer present in \eqref{Y1}. Taking the log-bilinear association model \eqref{Y1} instead of \eqref{Y2} weakens the distributional assumption while still including
the regression parameter $\beta$ up to a positive constant $\phi^{-1}$. In particular a linear hypothesis $Q\b = 0$ with a given matrix $Q$ is equivalent to $Q\t = 0$, and for a vector $c$ a one-sided hypothesis $c^T\b \leq 0$ is equivalent to $c^T \t \leq 0$.
\par\bigskip\noindent
A closer look at the relationship between generalized linear models and log-bilinear odds ratio models is given in section 6.2.
\par\bigskip\noindent
\textbf{Example 2:} Log-linear models for contingency tables\par \medskip\noindent
An important special case of example 1 are log-linear models for for contingency tables. 
If $X$ and $Y$ have finite support $\O_X=\{x_0,\;\ldots,\;x_J\}$ and $\O_Y = \{y_0,\;\ldots,\;y_K\}$ say, then the log-bilinear association model \eqref{C66} can be written as
\begin{equation}\label{1}
\psi_{jk}(\t)= \tilde{x}^T_j\t \tilde{y}_k,\hspace{1cm}\text{with}\hspace{1cm} \tilde{x}_j=h_X(x_{j}),\hspace{1cm} \tilde{y}_k=h_Y(y_{k}),
\end{equation}
or in matrix notation
\begin{equation}\label{2}
\psi (\t)= \tilde{X} \t \tilde{Y}^T\in\mathbb{R}^{J\times K},\hspace{1cm} \tilde{X} = (\tilde{x}_{jl}) \in\mathbb{R}^{J\times L_X},\hspace{1cm} \tilde{Y}=(\tilde{y}_{ki})\in\mathbb{R}^{K\times L_Y}.
\end{equation}
Then \eqref{Y5} reduces to a log-linear model for the probabilities $p_{jk}=p(x_j,\;y_k)$ namely
\begin{equation}\label{(1)}
\log p_{jk} = \a+ \r_j + \g_k + \tilde{x}^T_j \t \tilde{y}_k.
\end{equation}
with $\r_0 = \g_0 = 0$. \par\bigskip\noindent
Using Kronecker products the model \eqref{1} resp.\eqref{2} can be written as
\begin{equation}\label{3}
\begin{split}
\psi_{jk}(\t) = z^T_{jk}\vec{\t}\hspace{1cm}&\text{with}\hspace{1cm}z_{jk} = \tilde{y}_k\otimes \tilde{x}_j \in\mathbb{R}^L,\hspace{1cm}L=L_X L_Y\hspace{1cm}\text{resp.}\\
\vec{\psi}(\t) = Z\vec{\t}\hspace{1cm}&\text{with}\hspace{1cm}Z = \tilde{Y}\otimes \tilde{X} \in\mathbb{R}^{I\times L},\hspace{1cm}I=(J+1)(K+1)
\end{split}
\end{equation}
Note that the "interaction covariate" $z_{jk}$ is the vector representation of the $L_Y\times L_X$ matrix $\tilde{y}_k\tilde{x}_j^T$. The parameter $\t$ will be identifiable if and only if $\tilde{X}$ has rank $L_X$ and $\tilde{Y}$ has rank $L_Y$, i.e. $Z$ has rank $L$, and this will always be assumed.\par\bigskip\noindent
The saturated log-linear model imposes no restriction on the probabilities $p_{jk}$ and may be written as
\begin{equation}\label{4}
\log p_{jk} = \a + \r_j + \g_k + \psi_{jk}
\end{equation}
with constraints $\psi_{0k}=\psi_{j0}=0$. The model \eqref{(1)} can also be obtained by restricting the log-odds ratio table $\psi\text{°}=(\psi_{jk})_{j,k>0}$ to a linear subspace $\mathscr{Q}$ of $\mathbb{R}^{J\times K}$, namely $\mathscr{Q} = \{\tilde{X}\t \tilde{Y}^T | \t\in \mathbb{R}^{L_X\times L_Y}\}$. Hence log-bilinear association models are log-linear models where $\psi$ is restricted to a linear space, but the parameters $\r_1,\;\ldots,\;\r_J$ and $\g_1,\;\ldots,\;\g_K$ are not restricted (in order to leave the marginal distributions of $X$ and $Y$ unconstrained). 

\par\bigskip\noindent
\textbf{Example 3:} Multivariate linear logistic regression \par\medskip\noindent
Extending univariate logistic regression to the multivariate case, suppose $Y$ takes values in $\Omega_Y = \{0,\ 1,\;\ldots,\;K\}$, $K> 1$. Then $\mathscr{L}(Y | X = x)$ is a multinomial distribution $M_{K+1}(1,\pi(x))$ with $K +1$ classes and probabilities $\pi_k(x) = P(Y = k | X = x)>0$. Using the multivariate logistic transformation $\text{logit}\;\pi_k(x) = \log(\pi_k(x)/\pi_0(x))$, the multivariate linear logistic regression model is given by
\begin{equation}\label{(2.14)}
\text{logit}\;\pi_k(x) = \g_k +\tilde{x}^T \t_k,\hspace{1cm} k = 1,\;\ldots,\;K,
\end{equation}
where $\tilde{x}\in\RhX$ is as above a vector of formal covariates and $\g_k\in\mathbb{R}$, $\t_k\in\RhX$ are unknown parameters. Choosing $y_0= 0$, the log-odds ratio function is
\begin{equation}\label{Y3}
\psi(x, k) = \tilde{x}^T \t_k = \tilde{x}^T \t h_Y (k)=(h_Y (k)\otimes \tilde{x})^T \vec{\t},
\end{equation}
where $\t = (\t_1,\;\ldots,\;\t_K)$ is an $L_X\times K$ parameter matrix, and the function $h_Y :\Omega_Y \rightarrow \mathbb{R}^K$ maps $k > 0$ to the $k$th unit vector $e_k $ and $h_Y (0) = 0$. The model \eqref{(2.14)} is in fact equivalent to the log-bilinear association model \eqref{Y3} provided the parameters $\t_1,\;\ldots,\;\t_K$ are not restricted (cf.\ Osius 2004 \cite[sec. 4.2]{Met}). 
\par\bigskip\noindent

\noindent
\textbf{Example 4:} Multivariate linear regression\par\medskip\noindent 
Let $Y$ and $X$ be random vectors and suppose that the conditional distribution
of $Y$ given $X$ is multivariate normal,
\begin{equation}\label{(2.17)}
\mathscr{L}(Y | X = x) = N_{M_Y}(\mu_Y (x),\Sigma),
\end{equation}
such that the conditional covariance matrix $\Sigma$ is nonsingular and does not depend on $x$. From the conditional log-density
\begin{equation}
\log p(y | X = x)=-\frac{1}{2}\left[\log[(2\pi)^{M_Y} \text{det}(\Sigma)] + [y-\mu_Y (x)]^T \Sigma^{-1}[y-\mu_Y (x)]\right]
\end{equation}
the log-odds ratio function is 
\begin{equation}\label{(2.19)}
\psi(x, y) = [\mu_Y (x)-\mu_Y (x_0)]^T \Sigma^{-1}y. 
\end{equation}
The multivariate linear regression model
\begin{equation}\label{RM}
\mu_Y (x) = \alpha + \beta^T \tilde{x}
\end{equation}
with covariates $\tilde{x}$ and $L_X\times L_Y$ parameter matrix $\beta$ has a log-bilinear association
\begin{equation}
\psi(x, y) = \tilde{x}^T \t y
\end{equation}
with parameter matrix $\t = \beta \Sigma^{-1}$. The conditional covariance matrix $\Sigma$---and hence the parameter $\t$---may be recovered from the regression parameter $\beta$ and the (marginal) covariance matrices of $\tilde{X}$ and $Y$
\begin{equation}\label{(2.26)neu}
\Sigma = Cov(Y) - \beta^T Cov(\tilde{X})\beta,\hspace{1cm}    \t = \beta[Cov(Y) - \beta^T Cov(\tilde{X})\beta]^{-1}.
\end{equation}
Note that a linear hypothesis $C\b = 0$ is equivalent to the corresponding hypothesis $C\t = 0$, and the latter may be tested using the semiparametric association model \eqref{(2.19)} instead of the regression model \eqref{RM} with the additional distributional assumption \eqref{(2.17)}.

\section{Estimation}
We only give a brief overview of the estimation, for details see Osius, 2009 \cite[ch. 4]{Anals}. For a given data set $(x_i,y_i)$ with $i = 1,\;\ldots,\;n$ we want to estimate the association parameter $\t$ of the model \eqref{C6} under unconditional sampling from the joint distribution of $(X,Y)$ and conditional sampling of $Y$ given $X$ or vice versa. Not surprisingly the maximum likelihood estimator $\hat{\t}$ under any of these three sampling schemes may be obtained as a solution of the same estimating equation.

\subsection{Unconditional Sampling}
For unconditional sampling the data set $(x_i,y_i)$ is an independent sample from the joint distribution of $(X,Y)$. Suppose there are
$J+1>1$ different $x$-values and $K+1>1$ different $y$-values observed and denote the corresponding subsets of $\RX$ and $\RY$ by $\Omega_X^*=\{x_{(0)},\;\ldots,\;x_{(J)}\}$ and $\Omega_Y^*=\{y_{(0)},\;\ldots,\;y_{(K)}\}$. If $r_{jk}$ is the observed frequency of $(x_{(j)},y_{(k)})$, then the likelihood is
\begin{equation}
L_{XY} =\prod^J_{j=0} \prod^K_{k=0} p(x_{(j)},y_{(k)})^{r_{jk}} = L_{Y|X}\cdot L_X
\end{equation}
with a conditional and a marginal likelihood
\begin{equation}
L_{Y|X} =\prod^K_{k=0} \prod^J_{j=0} p(y_{(k)} | X = x_{(j)})^{r_{jk}},\hspace{1cm} L_{X} =\prod^J_{j=0} p_X (x_{(j)})^{r_{j+}}
\end{equation}
where the subscript \textquotedblleft$+$\textquotedblright\ indicates summation over the replaced index. The model does not restrict the marginal distributions of $X$ and $Y$ and hence the empirical densities with respect to counting measure,
\begin{align}
\hat{p}^Y(y_{(k)})&=r_{+k}/n \hspace{1cm} \text{for } k = 0,\;\ldots,\;K\label{Y4}\\
\hat{p}^X(x_{(j)})&=r_{j+}/n \hspace{1cm} \text{for } j = 0,\;\ldots,\;J
\end{align}
are the usual nonparametric estimators.
\par\bigskip\noindent
Interchanging $X$ and $Y$, we split the likelihood as $L_{XY} = L_{X|Y}\cdot L_Y$. Restricting $P^X$ and $P^Y$ to measures with finite support $\O^*_X$ and $\O^*_Y$ the likelihood $L_{XY}$ is a multinomial likelihood for the observed $(J+1)\times(K+1)$-contingency table $(r_{jk})$. And estimation of $\t$ is reduced to a multinomial model whose probabilities $p_{jk} = p(x_{(j)}, y_{(k)})$ satisfy the log-odds ratio model
\begin{equation}
\log\frac{p_{jk}p_{00}}{p_{j0} p_{0k}} = \psi_{\t}(x_{(j)}, y_{(k)})=\psi_{jk}(\t)\hspace{1cm} \text{for all $j$ and $k$}
\end{equation}
with respect to the reference values $x_0= x_{(0)}$ and $y_0= y_{(0)}$. The parametrization \eqref{Y5} now involves only a finite number of parameters
\begin{equation}\label{Y8}
\log p_{jk} = \rho_j + \g_k +\psi_{jk}(\t)-\log \left(\sum_j \sum_k \exp[\rho_j + \g_k +\psi_{jk}(\t)]\right),
\end{equation}
namely $\rho_j = \rho(x_{(j)})$, $\g_k = \g(y_{(k)})$ and $\t$ with $\rho_0 = \g_0 = 0$. Instead of maximizing $L_{XY}$, it is typically preferable to maximize either $L_{Y|X}$ or $L_{X|Y}$ using the parametrization of the conditional probabilities $p_{k|j} = p_{jk}/p_{j+}$ or $p_{j|k} = p_{jk}/p_{+k}$ given by \eqref{Y6} and \eqref{Y7}
\begin{equation}
\log p_{k|j} = \g_k +\psi_{jk}(\t)- \delta_j ,\hspace{1cm} \log p_{j|k} = \rho_j +\psi_{jk}(\t)-\varepsilon_k,
\end{equation}
where the parameters $\delta_j$, respectively $\varepsilon_k$, are determined by the remaining ones.

\subsection{Conditional Sampling} \label{3.2}
When sampling is conditional on values for $Y$ taken from $\Omega_Y^*= \{y_{(0)},\;\ldots,\;y_{(K)}\}$, say, then the data set $(x_i, y_i)$ with $i=1,\;\ldots,\;n$ is partitioned into $K +1$ independent subsamples given by the values of $y_i$, such that each subsample $(x_i)$ with $y_i=y_{(k)}$ is an independent sample from the conditional distribution $\mathscr{L}(X|Y=y_{(k)})$. Instead of maximizing the appropriate likelihood $L_{X|Y}$ we can equivalently maximize the unconditional likelihood $L_{XY}$ or even the \textquotedblleft reverse\textquotedblright\ conditional likelihood $L_{Y|X}$. The latter is preferable from a computational point of view, when the nuisance parameters $\g_k$ are less than those of $L_{X|Y}$, that is, for $K<L$. A dual argument applies if sampling is conditional on values for $X$ taken from $\Omega_X^*=\{x_{(0)},\;\ldots,\;x_{(J)}\}$.

\subsection{Log-bilinear Association}
In the log-bilinear association model \eqref{C6}, the odds ratios may be written as $\psi_{jk}(\t)= \tilde{x}^T_j\t \tilde{y}_k$ with $\tilde{x}_j = h_X(x_{(j)})$, $\tilde{y}_k = h_Y(y_{(k)})$ and a parameter matrix $\t\in\mathbb{R}^{L_X\times L_Y}$ or in matrix notation
\begin{equation}
\psi(\t) = \tilde{X}\t \tilde{Y}^T \in\mathbb{R}^{J\times K},\hspace{1cm} \tilde{X} = (\tilde{x}_{jl}) \in\mathbb{R}^{J\times L_X},\hspace{1cm} \tilde{Y} = (\tilde{y}_{kl}) \in\mathbb{R}^{K\times L_Y}.
\end{equation}
Then \eqref{Y8} reduces to a log-linear model for the probabilities $p_{jk}$,
\begin{equation}\label{(3.9)}
\log p_{jk} = \alpha +\rho_j +\g_k +\tilde{x}^T_j \t \tilde{y}_k
\end{equation}
induced by the covariates $\tilde{x}_j$, $\tilde{y}_k$. Hence results by Haberman, 1974 \cite[Ch. 2]{Hab} on the existence and uniqueness of maximum likelihood estimators in log-linear models apply. In particular the estimator $\hat{p} = (\hat{p}_{jk})$ is unique (if it exists) and the estimator $\hat{\t}$ is unique too, provided the parameter $\t$ is identifiable in the log-linear model \eqref{(3.9)}. As already noted in example 2, identifiability is equivalent to the conditions
\begin{equation}
\begin{split}
&\text{The $L_Y\times K$-matrix $\tilde{Y}^T = (\tilde{y}_1,\;\ldots,\;\tilde{y}_K)$ has rank $L_Y$ and}\\
&\text{the $L_X\times J$-matrix $\tilde{X}^T = (\tilde{x}_1,\;\ldots,\;\tilde{x}_J)$ has rank $L_X$}.
\end{split}
\end{equation}
This condition will be assumed here throughout. It will be satisfied if the sample is large enough, provided the functions $h_X$ and $h_Y$---and under conditional sampling the values $x_{(j)}$ resp. $y_{(k)}$---are properly chosen.\par\bigskip\noindent

\subsection{Log-linear Models for Contingency Tables}\label{(3.4)}
Since estimation of $\t$ in a log-bilinear association model can be reduced to estimation in a log-linear model we now have a closer look at the latter and continue with example 2.\par\bigskip\noindent
We now assume that $X$ and $Y$ have finite support $\O_X=\{x_0,\;\ldots,\;x_J\}$ resp. $\O_Y=\{y_0,\;\ldots,\;y_K\}$ and consider the usual sampling schemes for a $(J+1)\times(K+1)$ contingency table $R=(R_{jk})$ of random counts. 
The expected table will be denoted by $\mu=(\mu_{jk})=\mathbb{E}(R)$. It is important here that in all four sampling schemes the $I\times I$ covariance matrix $Cov(\vec{R})$ with $I=(J+1)(K+1)$ can be represented in terms of $D$-orthogonal projections onto a suitable linear subspace (cf.\ appendix A) where $D=diag\{\vec{\mu}\}$ is the diagonal matrix with diagonal $\vec{\mu}$. Furthermore the unit vector that stems from the $(J +1)\times(K+1)$ table having a one in the $(j,k)$th position and zeros otherwise will be denoted by $\vec{e}_{jk}$.
\par\bigskip\noindent

\subsubsection*{Multinomial Sampling}
Here we take an independent sample $(X_1,\;Y_1),\;\ldots,\;(X_n,\;Y_n )$ of size $n$ from the joint distribution of $(X, Y)$ and the $(J+1)\times(K+1)$-table $R = (R_{jk})$ of counts
\begin{equation*}
R_{jk} = \#\{i = 1,\;\ldots,\;n \;|\; X_i = x_j ,\; Y_i= y_k\}
\end{equation*}
follows a multinomial distribution
\begin{equation}
\mathscr{L}(\vec{R}) = M_{(J+1)(K+1)}(n, \vec{p})\hspace{1cm} \text{with}\hspace{1cm} p = (p_{jk})\tag{M}
\end{equation}
with $\mu_{jk}= n\cdot p_{jk}$. Define $\mathscr{D}=\ span\{\vec{e}_{++}\}$ \label{Q43}as the diagonal space that consists of all constant vectors in $\mathbb{R}^{I}$ and let $P_{\mathscr{D}}^D$ be the $D$-orthogonal projection onto the space $\mathscr{D}$, then (cf. Franke, 2010 \cite[sec. 2.2]{AF}; Habermann, 1974 \cite[(1.54)]{Hab}) 
\begin{equation}
Cov(\vec{R})=\ D - n^{-1} \vec{\mu}\vec{\mu}^T=\ D(\mathbb{I}-P_{\mathscr{D}}^D).
\end{equation}
The model \eqref{(1)} may also be written as a log-linear model for the expectations $\mu_{jk}$
\begin{equation}\label{(2)}
\log \mu_{jk}=\a' + \r_j + \g_k + \tilde{x}_j^T \t \tilde{y}_k 
\end{equation}
with $\a' = \a + \log n$.

\subsubsection*{Poisson Sampling}
Consider now an independent sample $(X_1,\;Y_1),\;\ldots,\;(X_N,\;Y_N)$ from the joint distribution of $(X,\;Y)$ where the sample size $N$ is an independent random variable having a Poisson distribution $Pois(\nu)$ with expectation $\nu$. Then the counts
\begin{equation*}
R_{jk} = \#\{i = 1,\;\ldots,\;N \;|\; X_i = x_j ,\; Y_i= y_k\}
\end{equation*}
are independent each having a Poisson distribution $Pois(\mu_{jk})$ with $\mu_{jk} = \nu p_{jk}$ and total expectation $\mu_{++}=\nu$. Hence the vector $\vec{R}$ has a product-Poisson distribution and we get the Poisson model
\begin{equation}
\mathscr{L}(\vec{R}) = \prod^J_{j=0}\prod^K_{k=0}\; Pois(\mu_{jk})\tag{P}
\end{equation}
with $p_{jk}=\ \;\mu_{jk}/\mu_{++}$ and $Cov(\vec{R})=D$. The model \eqref{(1)} may again be written as in \eqref{(2)} with $\a ' = \a +\log\l$.

\subsubsection*{Product Multinomial Sampling for Rows}
We now look at sampling conditional on $X$ where for each $j = 0,\;\ldots,\;J$ independent samples $X_{j1},\;\ldots,\;X_{jn_j}$ of size $n_j$ are taken from the conditional distribution $\mathscr{L}(Y|\;X = x_j)$. The rows $R_{j\cdot} = (R_{j0},\;\ldots,\;R_{jK})$ of the counts
\begin{equation*}
R_{jk} = \#\{i = 1,\;\ldots,\;n_j \;|\; Y_i= y_k\}
\end{equation*}
are independent for $j = 0,\;\ldots,\;J$ each with a multinomial distribution 
\begin{equation*}
\mathscr{L}(R_{j\cdot})=M_{K+1}(n_j,\; p^{|X}_j)\hspace{1cm}\text{with}\hspace{1cm} p^{|X}_{jk}= P\{Y=y_k|\; X=x_j\}.
\end{equation*}
Hence the vector $\vec{R}$ has a product-multinomial distribution and we get the product multinomial sampling for rows
\begin{equation}
\mathscr{L}(\vec{R})=\prod^J_{j=0} M_{K+1}(n_j,\; p^{|X}_j)\tag{MR},
\end{equation}
with $\mu_{jk}=n_jp^{|X}_{jk}$ and $Cov(\vec{R})=diag \left\{(\Sigma_j)_{j=0,\;\ldots,\;J}\right\}$ is a $(I\times I)$ block-diagonal matrix with blocks $\Sigma_j=Cov(R_{j\cdot})=diag\{\mu_{j\cdot}\} - \mu_{j+}^{-1}\mu_{j\cdot}\mu_{j\cdot}^T$ and $\mu_{j\cdot}$ the $j$th row of $\mu$. The columns of the $I\times(J+1)$ matrix $F=(\vec{e}_{0+},\;\ldots,\;\vec{e}_{J+})$ \label{Q42}span the row space $\mathscr{R}=span\{\vec{e}_{0+},\ldots , \vec{e}_{J+}\}$\label{Q44} which consists of all vectors arising from $(J+1)\times(K+1)$ tables with constant rows. Since $\left\langle \vec{e}_{j+},\vec{e}_{l+}\right\rangle_D=\delta_{lj}\cdot n_j$ (using Kronecker's $\delta$) the vectors $\vec{e}_{0+},\;\ldots,\;\vec{e}_{J+}$ are pairwise $D$-orthogonal. The covariance matrix of $\vec{R}$ can also be represented as (cf.\ Franke, 2010 \cite[sec. 2.3]{AF}; Habermann, 1974 \cite[(1.54)]{Hab}),
\begin{align}\label{(3.14)}
Cov(\vec{R})=diag \left\{(\Sigma_j)_j\right\}=\ \;&diag\{(diag\{\mu_{j\cdot}\} - \mu_{j+}^{-1}\mu_{j\cdot}\mu_{j\cdot}^T)_j\}=D(\mathbb{I}-P_{\mathscr{R}}^D).
\end{align}
Again the model \eqref{(1)} may be written in terms of the expectations as
\begin{equation}
\log \mu_{jk} = \a + \r_j' + \g_k + \tilde{x}^T_j\t \tilde{y}_k 
\end{equation}
with $\r_j'=\r_j+\log(n_j/p_{j+})$.

\subsubsection*{Product Multinomial Sampling for Columns}
Let us finally consider sampling conditional on $Y$ where for each $k = 0,\;\ldots,\;K$ we take independent samples $X_{k1},\;\ldots,\;X_{km_k}$ of size $m_k$ from the conditional distribution $\mathscr{L}(X|\; Y = y_k)$. The columns $R_{\cdot k}= (R_{0k},\;\ldots,\;R_{Jk})$ of the counts $R_{jk} = \#\{i = 1,\;\ldots,\;m_k \;|\; X_i= x_j\}$ are independent for $k = 0,\;\ldots,\;K$ each with a multinomial distribution. Hence the vector $\vec{R}$ has a product-multinomial distribution and satisfies the product multinomial sampling for columns
\begin{equation}
\mathscr{L}(\vec{R})=\prod^K_{k=0} M_{J+1}(m_k,\; p^{|Y}_k)\hspace{1cm}\text{with}\hspace{1cm} p^{|Y}_{kj}= P\{X=x_j|\; Y=y_k\}=p_{jk}/p_{+k}\tag{MC},
\end{equation}
$\mu_{jk}=m_k p^{|Y}_{kj}$ and $Cov(\vec{R})=diag \left\{(diag\{\mu_{\cdot k}\} - \mu_{+k}^{-1}\mu_{\cdot k}\mu_{\cdot k}^T)_{k=0,\;\ldots,\;K}\right\}$ 
with $\mu_{\cdot k}$ the $k$th column of $\mu$. 
The columns of the $I\times (K+1)$ matrix $G=(\vec{e}_{+0},\;\ldots,\;\vec{e}_{+K})$ span the column space $\mathscr{C}=span\{\vec{e}_{+0},\;\ldots,\;\vec{e}_{+K}\}$\label{Q46} which consists of all vectors arising from $(J+1)\times(K+1)$ tables with constant columns. Interchanging rows with columns, i.e.\ looking at the transposed table $R^T$, leads us back to the product model for rows and \eqref{(3.14)} yields
\begin{equation}
\begin {split}
Cov(\vec{R})=\ D(\mathbb{I}-P_{\mathscr{C}}^D).
\end{split}
\end{equation}

\subsubsection{Log-linear Models for the Expected Table}
In all four sampling schemes above, the expected table $\mu = \mathbb{E}(R)$ satisfies a log-linear model 
\begin{equation}\label{(5)}
\eta_{jk}= \a + \r_j + \g_k + \tilde{x}_j^T \t \tilde{y}_k\hspace{1cm}\text{respectively}\hspace{1cm}\vec{\eta}=\log\vec{\mu}\in\mathscr{H}
\end{equation}
with $\mathscr{H}$ a linear subspace of $\mathbb{R}^{I}$. 
\par\bigskip\noindent
Viewing $\tilde{x}_1,\;\ldots,\;\tilde{x}_J$ and $\tilde{y}_1,\;\ldots,\;\tilde{y}_K$ as "scores" assigned to the rows resp. columns, the above model appears as a generalization of the linear-by-linear association model in Agresti, 1990 \cite[sec. 8.1.1]{Agr} with vector-values scores instead of scalars. The above model may be rewritten as
\begin{align}
\eta_{jk} &= \a + \r_j + \g_k + z_{jk}^T\vec{\t} \hspace{1cm}\text{with}\label{5}\\
z_{jk} &= \tilde{y}_k\otimes\tilde{x}_j.\label{5a}
\end{align}
The vector $z_{jk}$ of dimension $L$ may be interpreted as an "interaction covariate" associated to $(j,k)$th cell of the $(J+1)\times(K+1)$-table and satisfies the constraints $\vec{z}_{j0}=\vec{z}_{0k}=0$. Although any log-linear model is of the form \eqref{5} it will only represent a log-bilinear association in our sense if the "covariate" $z_{jk}$ has a decomposition \eqref{5a}, which guarantees that $z_{jk}$ does not contain any information about the  association of $X$ and $Y$.
\par\bigskip\noindent
In the Poisson model (P) the parameters $\a'$, $\r_j$ and $\g_k$ are not restricted (cf.\ example 2) or equivalently, the marginal space $\mathscr{T}=\mathscr{R}+\mathscr{C}=\ span\{\vec{e}_{0+},\;\ldots,\;\vec{e}_{J+},\ \vec{e}_{+0},\;\ldots,\;\vec{e}_{+K}\}$ is a linear subspace of $\mathscr{H}$, and this will be assumed from now on. 
\par\bigskip\noindent
Given an observed table $r$ of counts the maximum likelihood estimator (in any of the four sampling schemes) $\vec{\hat{\mu}}=\hat{\mu}(\vec{r})\in\mathscr{M}=\exp[\mathscr{H}]$ of $\mu$ 
is the unique solution (provided there is one) of the same normal equation
\begin{equation}\label{C18}
P_{\mathscr{H}}\vec{\hat{\mu}}=P_{\mathscr{H}}\vec{r},
\end{equation}
cf.\ Haberman, 1974 \cite[ch. 2]{Hab} who also gives criteria for the existence of the estimate. In particular $\mathscr{T}\subset \mathscr{H}$ implies that $\hat{\mu}$ and $r$ have the same row and column totals 
\begin{align}
\hat{\mu}_{j+}&=\ r_{j+} &\text{for }j=0,\;\ldots,\;J,\\
\hat{\mu}_{+k}&=\ r_{+k} &\text{for }k=0,\;\ldots,\;K.\notag
\end{align}
The odds ratio parameter $\t$ is a function of $\eta$ resp. $\mu$ and will be estimated as the corresponding function. 
Conversely, $\hat{\mu}$ is the unique table determined by the log-odds ratios $\hat{\psi}_{jk} = \tilde{x}_j^T \hat{\t} \tilde{y}_k$ and the totals $r_{j+}$ and $r_{+k}$ of the observed table for all $j$ and $k$ (cf.\ Plackett, 1974 \cite[sec. 3.4] {Pla}).

\section{Asymptotic Covariance Matrices}
In this section---which contains the main results of this paper---we derive different representations for the (estimated) asymptotic covariance matrix $\Sigma_{\hat{\t}}$ of the estimator $\hat{\t}$. Here we assume that $Y$ has finite support and show in section 6 how the general case with arbitrary support for $Y$ can be reduced to finite support. We first look at log-linear models for contingency tables (example 2) where $X$ has finite support too. Then we consider the multivariate linear logistic regression model (example 3) with arbitrary support for $X$. Although the asymptotic covariance matrices arise from suitable asymptotic assumptions---and are only applicable given these assumptions---their estimates can always be computed for a given sample. And---using matrix algebra only---we are going to show that the different estimates considered here all result in the same matrix.
\subsection{Log-linear Models for Contingency Tables}\label{4.1}
Continuing our discussion in \ref{(3.4)} we consider a log-linear model given by $\eta\in\mathscr{H}$ with $\mathscr{T}\subset \mathscr{H}$ and assume \textit{any} of the four distribution models (M), (P), (MR) or (MC). The asymptotic normality of the estimates $\vec{\hat{\mu}}$ and $\vec{\hat{\eta}}$ given by Haberman, 1974 \cite[Th 4.4]{Hab}---for an asymptotic approach with fixed cells (i.e. $J$ and $K$ are fixed) and (suitably) increasing expectations $\mu_{jk}$ in each cell $(j,k)$---imply that the asymptotic covariance matrices of $\vec{\hat{\mu}}$ and $\vec{\hat{\eta}}$ are given by 
\begin{align}
\Sigma_{\hat{\mu}}\ =&D[P^{D}_{\mathscr{H}}-P^{D}_{\mathscr{N}}],\hspace{1cm}
\Sigma_{\hat{\eta}}\ =[P^{D}_{\mathscr{H}}-P^{D}_{\mathscr{N}}]D^{-1}\\ 
\text{with }\hspace{1cm}
&\mathscr{N}=\left.\begin{cases} \mathscr{D} &\text{for the model (M)},\\
														\mathscr{R} &\text{for the model (MR)},\\
														\mathscr{C} &\text{for the model (MC)},\\
														\{0\}&\text{for the model (P)}\end{cases}\right\}\subset\mathscr{T} \hspace{1cm}\text{and} \hspace{1cm} D=diag\{\vec{\mu}\}.
\end{align}
In each of the four sampling schemes the projection $P^{D}_{\mathscr{N}}Y$ is fixed by design, e.g. the row sums $Y_{j+}$ in (MR), and the distribution of $Y$ may be obtained from the Poisson model (P) by conditioning upon  $P^{D}_{\mathscr{N}}Y=c$ for a suitable $c$. To derive the asymptotic covariance matrix $\Sigma_{\hat{\t}}$ of $\vec{\hat{\t}}$ we use the representation
\begin{equation}\label{100}
\eta_{jk} = \a + \r_j + \g_k + \psi_{jk}(\t)\hspace{1cm}\text{with}\hspace{1cm} \psi_{jk}(\t) = z^T_{jk}\vec{\t}.
\end{equation}
Although in a log-bilinear association model $z_{jk}$ is given by \eqref{3} we derive the following results in appendix \ref{Proof_Th3} without this restriction and consider the particular case \eqref{3} separately. For later purpose we consider the compound parameter $\l=(\g\text{°},\t)$ with $\g\text{°}=(\g_k)_{k>0}$. We define further the $(JK \times L)$ matrix $Z\text{°}=(z_{jk}^T)_{j,k>0}$, the $(I \times JK)$ matrix $C$ through the columns $c_{jk} =\ \vec{e}_{jk}+\vec{e}_{00}-\vec{e}_{j0}-\vec{e}_{0k},\;j,k>0$ and the $(I \times K)$ matrix $B$ through the columns $b_k=\vec{e}_{0k}-\vec{e}_{00}$.
\par\bigskip\noindent
\fbox{ \begin{minipage}{14,2cm}
\begin{Th}\label{Th3}
In the log-linear model given by $\eta\in\mathscr{H}$ and \eqref{100} the asymptotic covariance matrix of the estimator $\vec{\hat{\l}}$ is given by
\begin{equation}\label{Q5}
\Sigma_{\hat{\l}}=\begin{pmatrix} B^T\\Z\text{°}^{-}C^T \end{pmatrix} \Sigma_{\hat{\eta}} \begin{pmatrix} B,&CZ\text{°}^{-T} \end{pmatrix}
\end{equation}
or in block notation
\begin{equation}\label{AA1}
\Sigma_{\hat{\l}}=\begin{pmatrix} \Sigma_{\hat{\g}\text{°}} & \Sigma_{\hat{\g}\text{°}\hat{\t}}\\
														 								 \Sigma_{\hat{\t}\hat{\g}\text{°}} & \Sigma_{\hat{\t}}	\end{pmatrix}.
\end{equation} 
In particular the asymptotic covariance matrix of $\vec{\hat{\t}}$ is given by
\begin{equation} \label{Sigma333}
\Sigma_{\hat{\t}}=Z\text{°}^{-}C^T P_{\mathscr{H}}^D D^{-1}CZ\text{°}^{-T}
\end{equation}
and does not depend on the space $\mathscr{N}$ (which determines the sampling scheme).
\end{Th}
\end{minipage} }\par\bigskip\noindent
\textit{Remark:} The above representation of $\Sigma_{\hat{\t}}$ contains in $D$ the vector $\vec{\mu}$ of expectations which depends on the sampling scheme. However the estimate $\vec{\hat{\mu}}$---and hence corresponding estimate $\hat{\Sigma}_{\hat{\t}}$ of $\Sigma_{\hat{\t}}$---is the same in the sampling schemes (P), (M), (MR) and (MC) and can be recovered from $\hat{\t}$ and the row and column totals of the observed table. 

\subsection{An Explicit Representation of $\Sigma_{\hat{\t}}$}\label{5.4}
To get a more explicit representation of the asymptotic covariance matrix $\Sigma_{\hat{\t}}$ in terms of the vectors $z_{jk}$ we first eliminate the projection $P_{\mathscr{H}}^D$ in \eqref{Sigma333} and obtain the representation (cf. appendix \ref{Proof_Th4}).\par\bigskip\noindent
\fbox{ \begin{minipage}{14,2cm}
\begin{Th}\label{Th4}
In the log-linear model given by \eqref{100} the asymptotic covariance matrix of the estimator $\vec{\hat{\t}}$ is
\begin{equation}\label{OR2}
\Sigma_{\hat{\t}}=\ \ (Z\text{°}^{T}(C^TD^{-1}C)^{-1}Z\text{°})^{-1}
\end{equation}
with $D=diag\{\vec{\mu}\}$ and $Z\text{°}=(z_{jk})_{j,k>0}$. The matrix $C^TD^{-1}C$ has for $j\neq l$ and $k\neq m$ the following elements\par\bigskip\noindent
\begin{equation}\label{AA4}
\begin{split}
(C^TD^{-1}C)_{jk,jk}&=\mu_{jk}^{-1}+ \mu_{00}^{-1} + \mu_{j0}^{-1} + \mu_{0k}^{-1}\\
(C^TD^{-1}C)_{jk,lm}&=\mu_{00}^{-1}\\
(C^TD^{-1}C)_{jk,jm}&=\mu_{00}^{-1} +\mu_{j0}^{-1}\\
(C^TD^{-1}C)_{jk,lk}&=\mu_{00}^{-1} +\mu_{0k}^{-1}.
\end{split}
\end{equation}
\end{Th}
\end{minipage}}\par\bigskip\noindent
The remark to theorem \ref{Th3} still applies here. This compact form of $\Sigma_{\hat{\t}}$ is helpful to evaluate the influence of the covariates and the estimates on the asymptotic covariance matrix $\Sigma_{\hat{\t}}$.
\par\bigskip\noindent
\textbf{Example (saturated model):} 
For the saturated model $\mathscr{H}=\mathbb{R}^{(J+1)\times(K+1)}$ the matrix $Z\text{°}$ is the identity matrix. Hence
\begin{equation}
\Sigma_{\hat{\t}}=C^T D^{-1} C.
\end{equation}
and its estimate $\hat{\Sigma}_{\hat{\t}}$ can be evaluated from \eqref{AA4} with $\mu$ replaced by the observed table $r$. 
In particular for a $2\times 2$ contingency table $R$ with $\Omega_X=\{0,1\}$ and $\Omega_Y=\{0,1\}$, i.e. $J=K=1$, we get the scalar $\hat{\Sigma}_{\hat{\t}}=r_{11}^{-1}+ r_{00}^{-1} + r_{10}^{-1} + r_{01}^{-1}$ which is well known as the asymptotic variance of the estimator of the log-odds ratio parameter $\t$.$\hspace*{9.3cm}\square$
\par\bigskip\noindent
In appendix \ref{Proof_Th5} we derive another representation of $DP_{\mathscr{T}^{\bot_D}}^D$ and hence of $\Sigma_{\hat{\t}}$, which will be used to prove 
\par\bigskip\noindent
\fbox{ \begin{minipage}{14,2cm}
\begin{Th}\label{Th5}
In the log-linear model given by \eqref{100} the asymptotic covariance matrix of the estimator $\vec{\hat{\t}}$ can be written in terms of the covariance matrix $Cov_{\text{\tiny{MR}}}(\vec{R})=DP_{\mathscr{R}^{\bot_D}}^D$ for the sampling scheme (MR), cf.\ \eqref{(3.14)}, $E=\ (\vec{e}_{+1},\;\ldots,\;\vec{e}_{+K})$ and the $(J+1)\times(K+1)$ matrix $Z=(z_{jk})$ as
\begin{equation}
\Sigma_{\hat{\t}}\ =\left[Z^T Cov_{\text{\tiny{MR}}}(\vec{R})Z-Z^TCov_{\text{\tiny{MR}}}(\vec{R})E(E^TCov_{\text{\tiny{MR}}}(\vec{R})E)^{-1}E^TCov_{\text{\tiny{MR}}}(\vec{R})Z\right]^{-1}
\end{equation}
for the sampling schemes (P), (M), (MR) and (MC).
\end{Th}
\end{minipage}}\par\bigskip\noindent
Again, the remark to theorem \ref{Th3} applies. This representation has been used to evaluate the covariance matrix for the special cases $K=1$ and $K=2$ (which also apply to linear logistic regression as remarked in \ref{4.4}), cf.\ Franke, 2010 \cite[sec. 5.1.3]{AF}.

\subsection{$\Sigma_{\hat{\t}}$ in Log-bilinear Association Models}
In the log-bilinear model \eqref{3} the matrix $Z$° is the Kronecker product of $\tilde{X}\text{°}=(\tilde{x}_{jl})_{j>0, l=1,...L_x}$ and $\tilde{Y}\text{°}=(y_{kl})_{k>0, l=1,...Ly}$, i.e. $Z\text{°}=\tilde{Y}\text{°}\otimes\tilde{X}\text{°}$. From the properties of Kronecker's product the left inverse $Z$°$^-$ can be obtained from the left inverses $\tilde{X}\text{°}^-$ and $\tilde{Y}\text{°}^-$ of $\tilde{X}\text{°}$ and $\tilde{Y}\text{°}$ as
\begin{equation}
Z\text{°}^- =\tilde{Y}\text{°}^-\otimes\tilde{X}\text{°}^-,\hspace{1cm} Z\text{°}^{-T} =\tilde{Y}\text{°}^{-T}\otimes\tilde{X}\text{°}^{-T}.
\end{equation}
Theorem \ref{Th3} applied to a log-bilinear association model gives
\begin{equation}\label{102}
\Sigma_{\hat{\t}}=(\tilde{Y}\text{°}^-\otimes\tilde{X}\text{°}^-)C^T P_{\mathscr{H}}^D D^{-1}C(\tilde{Y}\text{°}^{-T}\otimes\tilde{X}\text{°}^{-T})
\end{equation}
and theorem \ref{Th4} yields\par\bigskip\noindent
\fbox{ \begin{minipage}{14,2cm}
\begin{Co}\label{Th6}
In the log-bilinear model given by \eqref{3} the asymptotic covariance matrix of the estimator $\vec{\hat{\t}}$ is
\begin{equation}\label{AA5}
\begin{split}
\Sigma_{\hat{\t}}=((\tilde{Y}\text{°}^T\otimes\tilde{X}\text{°}^T)(C^TD^{-1}C)^{-1}(\tilde{Y}\text{°}\otimes\tilde{X}\text{°}))^{-1}
													 \end{split}
\end{equation}
with $D=diag\{\vec{\mu}\}$.
\end{Co}
\end{minipage}}\par\bigskip\noindent
The $(J+1)(K+1)\times L$ matrix $Z=(z_{jk})$ is the Kronecker product $Z=\tilde{Y}\otimes\tilde{X}$ and theorem \ref{Th5} leads to\par\bigskip\noindent
\fbox{ \begin{minipage}{14,2cm}
\begin{equation}
\begin{split}
&\Sigma_{\hat{\t}}=
\left[(\tilde{Y}^T\otimes\tilde{X}^T) \left(Cov_{\text{\tiny{MR}}}(\vec{R})-Cov_{\text{\tiny{MR}}}(\vec{R})E(E^TCov_{\text{\tiny{MR}}}(\vec{R})E)^{-1}E^TCov_{\text{\tiny{MR}}}(\vec{R})\right)(\tilde{Y}\otimes\tilde{X})\right]^{-1}\\
&\text{with}\hspace{0.5cm} Cov_{\text{\tiny{MR}}}(\vec{R})=DP_{\mathscr{R}^{\bot_D}}^D.
\end{split}
\end{equation}
\end{minipage}}\par\bigskip\noindent

\subsection{Multivariate Linear Logistic Regression with Sampling Conditional on~$X$}\label{4.4}
Consider now the more general case where $X$ is a random vector with arbitrary support, but $Y$ still having finite support $\Omega_Y=\{y_0,\ldots,y_K\}$. We assume a sampling scheme conditional on $X$ and choose $J+1$ different values $\Omega_X^*=\{x_0,\;\ldots,\;x_J\}$ of $X$. For each $j=0,\;\ldots,\;J$ an independent subsample $Y_{ji}$ with $i=1, ..., n_j$ is drawn from the conditional distribution of $Y$ given $X=x_j$ and as before the counts for $y_k$ in this subsample are denoted by $R_{jk} = \#\{i\;|\;Y_{ji}=y_k\}$. The resulting distribution model for the contingency table $R$ is the product multinomial sampling for rows with conditional probabilities $\pi_{jk|X}= P\{Y=y_k \;|\; X=x_j\}$ that are specified through the multivariate linear logistic regression model
\begin{equation}\label{1000}
\text{logit}_k(\pi_{j|X})=\ \g_k+  z_{jk}^T\vec{\t}\hspace{1cm}\text{for }j=0,\;\ldots,\;J,\ k=1,\;\ldots,\;K
\end{equation}
with arbitrary covariates $z_{jk}$ satisfying $z_{j0}=z_{0k}=0$ for all $j$ and all $k$. Note that the following statements not only hold for bilinear odds-ratio models with $z_{jk}=h_Y(k)\otimes \tilde{x}_j$ from \eqref{Y3} but also for the more general model \eqref{1000}. \par\bigskip\noindent
The log-likelihood with respect to $\l=(\g,\vec{\t})$ is
\begin{equation}
\log L_{Y|X}(\l)=\ \sum_{j=0}^J\sum_{k=0}^K R_{jk}\log \pi_{jk|X}
\end{equation}
and the score vector $U(\l)$ is its gradient
\begin{equation}
U(\l)=\ \text{D}_{\l} \log L_{Y|X}(\l)^T
\end{equation}
with covariance matrix given by the second derivative of $L_{Y|X}$
\begin{equation}\label{5.18}
Cov(U(\l)) = \mathbb{E}(-\text{D}^2_{\l\l}L_{Y|X}(\l)) = -\text{D}^2_{\l\l}L_{Y|X}(\l).
\end{equation}
It is well known that the inverse of this matrix is---under mild conditions---the asymptotic covariance matrix of the estimator $\hat{\l}$ when the total sample size increases. In appendix \ref{Proof_Th7} we prove a fundamental result that $Cov(U(\l))^{-1}$ coincides with the asymptotic covariance matrix $\Sigma_{\hat{\l}}$ for $\hat{\l}$ given in theorem \ref{Th3}, where $X$ had \textit{finite} support.\par\bigskip\noindent
\fbox{ \begin{minipage}{14,2cm}
\begin{Th}\label{Th7}
For sampling conditional on $X$ the inverse of the covariance matrix of the score vector $U(\l)$ is given by
\begin{equation*}
Cov(U(\l))^{-1}=\Sigma_{\hat{\l}}
\end{equation*}
with $\Sigma_{\hat{\l}}$ from theorem \ref{Th3}.
\end{Th}
\end{minipage} }\par\bigskip\noindent
Hence the estimate $\hat{\l}$ and its asymptotic covariance $\Sigma_{\hat{\l}}$ can be determined as if $X$ had finite support $\Omega_X^*$. In particular any statistical software package for multivariate linear logistic regression or log-linear models can be used to compute $\hat{\l}$ and the estimate $\hat{\Sigma}_{\hat{\l}}$ of $\Sigma_{\hat{\l}}$ as well as to perform further statistical analysis, like tests and confidence intervals. Furthermore the representations of the estimated asymptotic covariance matrix $\hat{\Sigma}_{\hat{\t}}$ given in \ref{4.1} and \ref{5.4} apply here too.\par\bigskip\noindent

\section{Arbitrary Support of $Y$ and Sampling conditional on $Y$}
So far we have assumed that $Y$ has finite support and we now consider the general case with arbitrary support for $Y$ and $X$. Although the maximum likelihood estimate $\hat{\t}$ of the association parameter $\t$ may be obtained by maximizing the likelihood for conditional or unconditional sampling, the stochastic properties of the latter depend on the sampling scheme. Let us consider sampling conditional on $Y$---which can be preferable from a practical point of view---and summarize properties of the estimate $\hat{\t}$, for details see Osius, 2009 \cite[sec. 5--7]{Anals}. It is convenient to represent the sample as a compound vector $X = (X_{ki})$ of independent random variables indexed by $k = 0,\;\ldots,\;K$ and $i = 1,\;\ldots,\;m_k$. Using the notations from \ref{3.2} without the parentheses in $y_{(k)}$ and $x_{(j)}$, each $X_{ki}$ is distributed as $X_k \sim \mathscr{L}(X | Y = y_k)$. Let $R_{jk} = \# \{ i\ |\ X_{ki}=x_j \}$ denote the frequency of $x_j$ in the subsample $(X_{ki})$. Then $R_{+k} = m_k$ is fixed and the empirical distribution on $\O_Y^* = \{y_0,\;\ldots,\;y_K\}$ is given by the proportions $\bar{m}_k = m_k/n$, where $n = m_+$ is the total sample size. Replacing in the joint distribution $P$ of $(X, Y)$ the marginal distribution of $Y$ by the empirical distribution \eqref{Y4} yields a joint distribution $P^*$ on $\RX \times \O^*_Y$ given by the density $p^*$ with respect to the product of $\nu_X$ and the counting measure on $\O^*_Y$:
\begin{equation*}
p^*(x, y_k)= \bar{m}_k \cdot p(x | Y = y_k)\hspace{1cm} \text{for all $x$, $k$.}
\end{equation*}
Denoting the conditional density of $X$ by
\begin{equation}
\begin{split}
p^*_k(x) =&\;p^*(y_k | X = x) = \frac{\bar{m}_k \cdot p(x | Y = y_k)}{p^{*X}(x)},
\end{split}
\end{equation}
equation \eqref{Y6} yields the parametrization $\log p^*_k (x) = \g^*_k +\psi_{\t} (x, y_k)-\delta^*(x)$ with nuisance parameters $\g^*_k = \g^*(y_k)$ and $\delta^*(x) = \log[\sum_l \exp(\g^*_l +\psi_{\t} (x, y_l))]$, hence
\begin{equation}\label{(5.2)}
p^*_k (x) = \frac{\exp(\g^*_k +\psi_{\t} (x, y_k))}{\sum_l \exp(\g^*_l +\psi_{\t} (x, y_l))}.
\end{equation}
From the constraints \eqref{constr2.3} we obtain $\g^*_0 = 0$, and the nuisance parameter is $\g^*= (\g^*_1,\;\ldots,\;\g^*_K )\in\mathbb{R}^K$. Finally, the logarithm of the conditional likelihood $L_{Y |X}$ may be written in terms of the compound parameter vector $\l = (\g^*,\vec{\t})\in\mathbb{R}^{K+L}$:
\begin{equation}\label{10}
\begin{split}
l(\l)=& \log L_{Y |X} =\sum^K_{k=0} \sum^{m_k}_{i=1}\log p^*_k(X_{ki}).
\end{split}
\end{equation}
The first and second derivative of $l(\l)$ are denoted by $\text{D}_{\l}l(\l)$ and $\text{D}_{\l\l}^2l(\l)$.
\par\bigskip\noindent
Let us briefly resume the asymptotic properties of the estimator $\hat{\l}=(\hat{\g}^*,\hat{\t})$. The asymptotic approach assumes that set $\Omega_Y^*= \{y_0,\;\ldots,\;y_K\}$ of conditional values will remain fixed while all subsample sizes $m_k$ tend to infinity with fixed ratios $\bar{m}_k = m_k/n > 0$ for all $n$ and $k$. Hence the nuisance parameter $\g^*$ and the conditional densities $p_k^*(x)=p^*(y_k|X=x)$ do not vary with $n$. The asymptotic unique existence of the estimator,the strong consistency of the sequence $\hat{\lambda}^{(n)}$ and its asymptotic normality can be derived under reasonable conditions. More precisely, using a block notation for the inverse of the information matrix $\mathbf{I}(\l)=-\mathbb{E}(\text{D}^2_{\l\l}l(\l))$, i.e.
\begin{equation}\label{(4.3a)}
\mathbf{I}(\l)^{-1}=\begin{pmatrix} (\mathbf{I}(\l)^{-1})_{\g\g} & (\mathbf{I}(\l)^{-1})_{\g\t}\\ (\mathbf{I}(\l)^{-1})_{\t\g} & (\mathbf{I}(\l)^{-1})_{\t\t} \end{pmatrix},
\end{equation}
the asymptotic distribution is given by (cf.\ Osius, 2009 \cite[Thm. 5]{Anals})
\begin{equation}
\begin{split}
\sqrt{n}[\vec{\hat{\t}}^{(n)}-\vec{\t}]\xrightarrow[n\rightarrow\infty]{\mathscr{L}}N(0,(\bar{\mathbf{I}}^{-1}(\l))_{\t\t})\hspace{1cm}\text{with}\hspace{1cm}
\bar{\mathbf{I}}(\l)=n^{-1}\mathbf{I}(\l).\label{Th2}
\end{split}
\end{equation}
The matrix
\begin{equation}\label{Q78}
\mathbf{J}(\lambda)=\ -\text{D}^2_{\l\l}l(\l)
\end{equation}
is a consistent estimator of $\mathbf{I}(\lambda)$ and hence
\begin{equation}\label{7.1}
\begin{split}
\vec{\hat{\t}}\;&{}_{\stackrel{{}^{\text{\normalsize{$\sim$}}}}{{}_{\text{as.}}}}\;N(\vec{\t},(\mathbf{J}^{-1}(\hat{\lambda}))_{\t\t})\hspace{1cm}\text{with}\\
(\mathbf{J}^{-1}(\lambda))_{\t\t}&=\left(\mathbf{J}(\lambda)_{\t\t}-\mathbf{J}(\lambda)_{\t\g}(\mathbf{J}(\lambda)_{\g\g})^{-1}\mathbf{J}(\lambda)_{\g\t}\right)^{-1}
\end{split}
\end{equation}
using a well known result for the inverse of a partitioned matrix: 
\begin{equation}\label{Matrix1}
A=\begin{pmatrix}  L & M \\ G & H \end{pmatrix} \quad\Leftrightarrow \quad
A^{-1}=\begin{pmatrix}  L^{-1}+L^{-1}MN^{-1}GL^{-1} & -L^{-1}MN^{-1} \\ -N^{-1}GL^{-1} & N^{-1} \end{pmatrix}
\end{equation}
with $N=H-GL^{-1}M$. 
\par\bigskip\noindent
Note that for an observed data set, the estimated covariance matrix $(\mathbf{J}^{-1}(\hat{\lambda}))_{\t\t}$---i.e\ $\lambda$ is replaced in \eqref{7.1} by $\hat{\lambda}$---is identical to the corresponding matrix under sampling conditional on $X$ (instead of $Y$). In this sense
the estimate $\hat{\t}$ and its estimated asymptotic normal distribution are invariant under sampling conditional on either $Y$ or $X$. Hence asymptotic inference (i.e.\ tests or confidence regions) for the association parameter $\t$ based on the asymptotic distribution \eqref{7.1} of the estimate $\hat{\t}$ is invariant under both conditional sampling schemes, too.
\par\bigskip\noindent
For an observed table $(r_{jk})$ the matrix $\mathbf{J}(\hat{\l})$ may be computed as if sampling had been conditional on $X$ (instead of $Y$). However, for sampling conditional on $X$ \eqref{5.18} and theorem \ref{Th7} imply
\begin{equation}
\mathbf{J}^{-1}(\hat{\l})=\hat{\mathbf{\Sigma}}_{\hat{\l}}
\end{equation}
with $\hat{\mathbf{\Sigma}}_{\hat{\l}}$ from the remark to theorem \ref{Th3}. In particular the estimated asymptotic covariance matrix of $\hat{\t}$ coincides with the estimate of \eqref{Sigma333}
\begin{equation}
(\mathbf{J}^{-1}(\hat{\l}))_{\t\t} = Z\text{°}^{-}C^T P_{\mathscr{H}}^{\hat{D}} \hat{D}^{-1}CZ\text{°}^{-T}\hspace{1cm}\text{with}\hspace{1cm}\hat{D}=diag\{\hat{\mu}\}.
\end{equation}
We note again, that the table $\hat{\mu}$ is uniquely determined by the row and column totals of the observed table $(r_{jk})$ and the estimate $\hat{\t}$.
\par\bigskip\noindent
Hence, the estimated asymptotic covariance matrix of $\hat{\t}$ for sampling conditional on $Y$ is the same as for the usual fixed cells asymptotics where $X$ and $Y$ had \textit{finite} support. And interchanging $X$ and $Y$ yields the same result for sampling conditional on $X$.

\section{Applications}
This section deals with some applications of our theoretical results. The covariance matrix $\Sigma_{\hat{\t}}$ (for which we have given several representations) is not only needed to analyze a given sample by means of log-bilinear association models but also to investigate the properties of such an analysis, mainly the power of the tests involved and the calculation of the necessary sample size to achieve sufficient power. 
We first address power and sample size issues for unconditional and conditional sampling. And finally we have a closer look at generalized linear models with canonical link (in particular linear and log-linear models) and discuss the advantage of using the more general log-bilinear odds ratio models instead

\subsection{Power and Sample Size Issues}\label{ch6.1}
Suppose we wish to test a linear hypothesis $H_0:\ Q\t = 0$ for a given matrix $Q$ against the alternative $H:\ Q\t \neq 0$ using the usual test based on the asymptotic normal distribution of $Q\hat{\t}$. As a typical example, suppose $X=(X',X'')$ consists of two blocks and we wish to test the hypothesis $H_0$ that $X''$ and $Y$ are independent, which is often of primary interest. Using separate functions $\tilde{x}'=h_X'(x')$ and $\tilde{x}''=h_X''(x'')$ such that $\tilde{x}=(\tilde{x}',\tilde{x}'')$ and the block notation $\t=(\t',\t'')$, the above hypothesis of independence is equivalent to the linear hypothesis $H_0:\t''=0$. If in addition $Y=(Y',Y'')$, a similar argument shows that the hypothesis ``$X''$ \textit{and} $Y''$ \textit{are independent}'' is a linear hypothesis too.
\par\bigskip\noindent
The asymptotic power of the test of $H_0:\ Q\t=0$ may be computed from the covariance matrix of the estimator $\hat{\t}$ using one of the above representations of $\Sigma_{\hat{\t}}$. We first look at contingency tables, i.e. both $X$ and $Y$ have finite support, and consider unconditional and conditional sampling separately.
\par\bigskip\noindent
\textbf{Unconditional (multinomial) sampling for contingency tables:} In the multinomial sampling (M) the vector expectations is given by $\vec{\mu}=n\vec{p}$ and from corollary \ref{Th6} we get
\begin{equation}\label{(6.1)}
\Sigma_{\hat{\t}}^{-1} = n (\tilde{Y}\text{°}^T\otimes\tilde{X}\text{°}^T)(C^T diag^{-1}\{\vec{p}\}C)^{-1}(\tilde{Y}\text{°}\otimes\tilde{X}\text{°})\hspace{1cm} 
\end{equation}
using \eqref{AA4} to evaluate $C^Tdiag^{-1}\{\vec{p}\}C$. The matrices $X\text{°}$ and $Y\text{°}$ contain only the known values $\tilde{x}_j$ and $\tilde{y}_k$, but the joint density $p$ additionally depends on $\t$. For a given value $\t'$ of interest from the alternative we wish to compute the asymptotic power of the test either retrospectively (i.e. after the sample has been drawn) or prospectively to obtain an optimal design for the study. Since $X\text{°}$ and $Y\text{°}$ are already known we only
have to find the joint density $p'$ corresponding to $\t'$ and the marginal probabilities $p_{j+}$ and $p_{+k}$. This unique $p'$ can be obtained by an iterative proportional fitting procedure (cf.\ Sinkhorn, 1967 \cite{Sink}). Alternatively, $p'$ can be found by fitting the log-linear model \eqref{(1)} under the constraint $\t=\t'$ to an "observed" table $r'$ with marginals $r_{j+}'= p_{j+}$ and $r_{+k}'=p_{+k}$, e.g. $r_{jk}'= p_{j+}p_{+k}$. Using $p'$ instead of $p$ in \eqref{(6.1)} yields
\begin{equation}\label{(6.3)}
\Sigma_{\hat{\t}}'^{-1} = n\;(\tilde{Y}\text{°}^T\otimes\tilde{X}\text{°}^T)(C^Tdiag^{-1}\{\vec{p}\;'\}C)^{-1}(\tilde{Y}\text{°}\otimes\tilde{X}\text{°}).
\end{equation}
from which the asymptotic power of the test can be obtained.$\hspace{5cm}\square$
\par\bigskip
\noindent
\textbf{Conditional sampling for contingency tables:} Sampling conditional on $X$ leads to product multinomial sampling for rows (MR) where the expectations are given by $\mu_{jk}' = n_j p_{jk}'^{|X}$ with conditional probability $p_{jk}'^{|X}=p_{jk}'/p_{j+}'$. Using the total sample size $n = n_+$ and the relative sample sizes $\bar{n}_j= n_j/n$---which are typically fixed in advance, e.g. $\bar{n}_j=n/(K+1)$ in a balanced design---we get $\vec{\mu}'=n\vec{p}'^*$ with $p'^*_{jk}= \bar{n}_j p'_{jk}/p'_{j+}$ and hence
\begin{align}\label{(6.4)}
\Sigma_{\hat{\t}}'^{*-1} =& n\;(\tilde{Y}\text{°}^T\otimes\tilde{X}\text{°}^T)(C^Tdiag^{-1}\{\vec{p}\;'^*\}C)^{-1}(\tilde{Y}\text{°}\otimes\tilde{X}\text{°}).
\end{align}
The density $p'^*$ arises from $p'$ by replacing the marginal distribution of $X$ with the empirical distribution of $X$ given by the proportions $\bar{n}_j$. Note however, that the marginal distribution of $Y$ changes when passing from $p'$ to $p'^*$, i.e. $p_{+k}'=p_{+k}$ differs from $p_{+k}'^*$. Consequently the joint distribution $p'^*$ is not determined by $\t'$, $\bar{n}_j$ and $p_{+k}$ (for all $j,k$) alone, but still depends on the marginal distribution of $X$ although sampling is conditional on $X$.
\par\bigskip\noindent
The matrix \eqref{(6.4)}---and hence the power of the test---depends not only on the total sample size, but also on the proportions $\bar{n}_j$ which may be chosen to maximize the power.
\par\bigskip\noindent
And for sampling conditional on $Y$, i.e. the model (MC), we get the same representation \eqref{(6.4)} with $n=m_+$, $\bar{m}_k=m_k/n$ and $p_{jk}'^*= \bar{m}_k p_{jk}'/p_{+k}'$ . Again the power of the test may be maximized with respect to the proportions $\bar{m}_k$. And if conditional sampling on $X$ or $Y$ are both possible then one can choose the sampling design with the highest power for the test.$\hspace{13.7cm}\square$
\par\bigskip\noindent
To determine the total sample size $n$ necessary to achieve a wanted power, we only have to increase $n$ in \eqref{(6.3)} resp. \eqref{(6.4)} until the given power is reached. The above consideration only apply when both $X$ and $Y$ have finite range. However the distributions of $X$ and $Y$ can always be approximated by distributions with finite support, e.g. by grouping or rounding. And using the discrete approximations to compute the power should be sufficiently accurate for practical purposes.

\subsection{Generalized Linear Models With Canonical Link vs. Log-Bilinear Odds Ratio Models}
In example 1 we have already seen that generalized linear models with canonical link function are log-bilinear odds ratio models. However the latter models do not assume that the conditional distributions belong to the exponential family \eqref{(2.9)}. We now explore in more detail the relationship between these regression and association models. Keeping the notation from section \ref{sec2} we suppose that $Y$ is \textit{univariate} with support $\O_Y\subset \mathbb{R}$. We consider the log-bilinear odds ratio model with respect to the identity map $h_Y = id$ on $\mathbb{R}$---e.g. $\tilde{y} = y$--- 
\begin{equation}\label{Verlobung}
\psi_{\t}(x,y) = \tilde{x}^T\t y\hspace{1cm} \text{for all } x,\ y.
 \end{equation}
This model does not restrict the marginal distributions $P^X$ and $P^Y $ of $X$ and $Y$. But we assume that $Cov(\tilde{X})$ is positive definite and $0< \sigma_Y^2=Var(Y) <\infty$ which guarantees for \textit{any} $\t\in\mathbb{R}^{L_X}$ the existence of a unique joint distribution $P$ with \eqref{Verlobung} and marginals $P^X$ and $P^Y$. The logarithm of the conditional density \eqref{Y6} of $Y$ given $X = x$ may now be written
\begin{equation}
\begin{split}\label{(6.5)}
\log p(y|x)&=\g(y) + \tau y - \kappa(\tau)\hspace{1cm} \text{with}\\
\tau&=\tilde{x}^T \t\hspace{1cm} \text{and}\\
\kappa(\tau)&=\log \int \exp(\g(y) + \tau y)d\nu_Y(y).
\end{split}
\end{equation}
Although this density looks like a member of an exponential family with canonical parameter $\tau$, it need not be a density for any value of $\tau$ other than $\tilde{x}^T \t$. However the expectation and variance of the conditional distribution are still given by the derivatives of $\kappa$
\begin{equation}\label{(6.66)}
\begin{split}
\mathbb{E}(Y|X=x)&=\kappa'(\tau)=\kappa'(\tilde{x}^T \t)=:\mu_x(\t),\\
Var(Y|X=x)&=\kappa''(\tau)=\kappa''(\tilde{x}^T \t)=:\sigma_x^2(\t)
\end{split}
\end{equation}
provided the following regularity condition holds which allows interchanging differentiation with integration
\begin{equation}\label{RC}
\begin{split}
\text{D}_\tau\int \exp(\g(y) +\tau y) d\nu_Y(y) &= \int\left[\text{D}_\tau\exp(\g(y) +\tau y)\right] d\nu_Y(y) ,\\
\text{D}^2_{\tau\tau}\int \exp(\g(y) +\tau y) d\nu_Y(y)& = \int\left[\text{D}^2_{\tau\tau} \exp(\g(y) +\tau y)\right] d\nu_Y(y).
\end{split}
\end{equation}
The derivative of the conditional expectation with respect to $\t$ is
\begin{equation}\label{(6.77)}
\mu'_x(\t)=\kappa''(\tilde{x}^T \t)\tilde{x}^T=\sigma_x^2(\t)\tilde{x}^T.
\end{equation}
We will now see how the linear resp. log-linear or logistic regression model emerges from the association model when the respective structure for the conditional variance is assumed.

\subsubsection{Linear Regression}
Now $Y$ has a \textit{continuous} distribution with support $\O_Y= \mathbb{R}$ and we assume that the conditional variance is constant and positive
\begin{equation}\label{(6.7)}
\sigma_x^2(\t)=\sigma^2>0\hspace{1cm} \text{for all $x$ and $\t$},
\end{equation}
which is a common assumption in linear regression models. Then the derivative $\mu_x'(\t)$ does not depend on $\t$ and hence the conditional expectation may be written as
\begin{align}
\mu_x(\t)&= \beta_0 + \beta^T \tilde{x}\hspace{1cm}\text{for all $x$ with}\label{(6.8)}\\
\beta&=\sigma^2\t\label{(6.9)}
\end{align}
and some constant $\beta_0\in\mathbb{R}$. Conversely, the linear model \eqref{(6.8)} and \eqref{(6.9)} together
imply \eqref{(6.7)} in view of \eqref{(6.77)}. From \eqref{(6.7)} and \eqref{(6.8)} one easily obtains
\begin{equation}\label{(6.10)}
\begin{split}
\mathbb{E}(Y)=&\beta_0 + \beta^T \mathbb{E}(\tilde{X})\\
\sigma^2=&\sigma_Y^2 - \beta^T Cov(\tilde{X}) \beta=\sigma_Y^2 - ||\beta||^2_{Cov(\tilde{X})}
\end{split}
\end{equation}
using the norm induced by $Cov(\tilde{X})$ (cf. appendix A). This in turn gives the odds ratio parameter $\t$ in terms of the regression parameter $\beta$ and the second order moments of the (marginal) distributions of $\tilde{X}$ and $Y$
\begin{equation}\label{(6.11)}
\t = [\sigma_Y^2 - ||\beta||^2_{Cov(\tilde{X})}]^{-1} \beta,
\end{equation}
which coincides with \eqref{(2.26)neu} for univariate $Y$. In order to recover the regression parameter $\beta$ from $\t$---given $\sigma_Y$ and $Cov(\tilde{X})$---we consider the norms
\begin{equation}\label{(6.13)}
||\t||_{Cov(\tilde{X})} = [\sigma_Y^2 - ||\beta||^2_{Cov(\tilde{X})}]^{-1} ||\beta||_{Cov(\tilde{X})}=f(||\beta||_{Cov(\tilde{X})}).
\end{equation}
The function $f(u) = u / (\sigma_Y^2 - u^2)$ defined for $u^2\neq \sigma_Y^2$ with derivative $f'( u)= (\sigma_Y^2+ u^2)/ (\sigma_Y^2 - u^2)^2$ is strictly increasing for $0 \leq u < \sigma_Y$ from $f(0) = 0$ to its left-sided limit $f(\sigma_Y -) =\infty$. Hence $f$ has an inverse $f^{-1}: [0,\infty)\;\longrightarrow\;[0,\sigma_Y)$ given by
\begin{equation}
f^{-1}(v)=\begin{cases} 0 & \text{for } v =0,\\
\frac{1}{2} v^{-1} \left[ \sqrt{1+4v^2\sigma_Y^2}\right] & \text{for } v> 0 .\end{cases}
\end{equation}
Now we obtain $||\beta||_{Cov(\tilde{X})}=f^{-1}(||\t||_{Cov(\tilde{X})})$ from \eqref{(6.13)}, which inserted in \eqref{(6.11)} yields $\beta$ in terms of $\t$ and the second order moments of the (marginal) distributions of $\tilde{X}$ and $Y$
\begin{equation}\label{(6.15)}
\beta = \left[\sigma_Y^2 - f^{-1}(||\t||_{Cov(\tilde{X})})^2\right] \t .
\end{equation}
From the above discussion the log-bilinear association model \eqref{Verlobung} appears as a generalization of the classical linear model which---in addition to \eqref{(6.7)}---assumes the conditional distribution to be normal
\begin{equation}
\mathscr{L}(Y|X=x)=N(\mu_x(\t),\sigma^2).
\end{equation}
Furthermore, the linear model only leaves the marginal distribution of $X$ unconstrained but introduces a connection between the marginal distributions of $X$ and $Y$, e.g.\ through \eqref{(6.10)}.
\par\bigskip\noindent
As already mentioned, using the association model \eqref{Verlobung} instead of the regression model \eqref{(6.8)} with \eqref{(6.7)} also allows asymptotic inference about $\beta$---for sampling conditional on either $X$ or $Y$---because $\t$ and $\beta$ only differ by the positive (unknown) constant $\sigma^2$.
\par\bigskip\noindent
Furthermore a one-sided hypothesis $H_0:\ c^T\beta\leq 0$ for a given vector $c$ is equivalent to $H_0:\ c^T\t\leq 0$ and a linear hypothesis $H_0:\;Q\beta = 0$ for a given matrix $Q$ is equivalent to $H_0:\;Q\t = 0$. To compute the power of the corresponding test for a given value $\beta'$ under the alternative we have to assume realistic values for the variance $\sigma_Y^2$ and the covariance matrix $Cov(\tilde{X})$ in order to get the corresponding values of $\sigma'^2$ and $\t'$ from \eqref{(6.10)} and \eqref{(6.9)} for $\beta=\beta'$. Then the considerations in section \ref{ch6.1} can be applied to obtain---for the intended sampling scheme---the corresponding covariance matrix $\Sigma'_{\hat{\t}}$ which allows the computation of the power and the necessary sample size to achieve a given power.
\par\bigskip\noindent
Using the log-bilinear odds ratio model determined by \eqref{Y7} and \eqref{C66}  has two advantages over the usual linear regression model given by \eqref{(6.7)} and \eqref{(6.8)}. First, no assumptions about the conditional distribution of $Y$ given $X$ are needed and in particular, \eqref{(6.7)} need not hold. And second, sampling may be conditional on $Y$ instead of $X$, which may be preferable from a practical point of view or to achieve a higher power. 
\par\bigskip\noindent
However, even if the linear model holds and if the marginal variance $\sigma_Y^2$ and the covariance matrix $Cov(\tilde{X})$ are known---or consistent estimates are available, e.g.\ from previous studies---then a plug-in estimator $\hat{\beta}$ of $\beta$ can be obtained from \eqref{(6.15)} and $\hat{\t}$. Furthermore the asymptotic normality of $\hat{\t}$ provides the asymptotic normal distribution of $\hat{\beta}$ by the delta-method.

\subsubsection{Log-Linear Regression}
We now consider the case where $Y$ is discrete with support $\O_Y = \mathbb{N} \cup \{0\}$ and assume
\begin{equation}\label{(PV)} 
\sigma_x^2(\t)=\mu_x(\t)>0 \hspace{1cm}\text{for all $x$ and $\t$,}
\end{equation}
i.e. the Poisson variance function applies. Then by \eqref{(6.66)} $\kappa'(\tilde{x}^T\t) = \kappa''(\tilde{x}^T\t)$ for all $x$ and $\t$---which in turn implies $\kappa'(\tilde{x}^T\t) = \exp(\beta_0 + \tilde{x}^T\t) + c$ for some constants $\beta_0,\ c\in\mathbb{R}$. If the expectation $\mu_x(\t)$ is allowed to take any positive value, then $c$ must be zero and we get the familiar log-linear model
\begin{equation} \label{(6.16)}
\log \mu_x(\t)=\beta_0 + \tilde{x}^T\beta\hspace{1cm} \text{with }\beta=\t.
\end{equation}
Hence the association model \eqref{Verlobung} appears as a generalization of the log- linear model which---in addition to \eqref{(PV)}---restricts the conditional distributions to Poisson distributions
\begin{equation}
\mathscr{L}(Y|X=x)=Pois(\mu_x(\t)).
\end{equation}
Since $\beta=\t$ asymptotic inference about the regression parameter $\beta$ of the log-linear model \eqref{(6.16)} may also be obtained from the more general association model which imposes no restriction on the conditional distribution of $Y$ given $X$, e.g. \eqref{(PV)}, and where sampling may be conditional on either $X$ or $Y$.

\subsubsection{Logistic Regression}
Looking finally at a binary random variable $Y$ with support $\O_Y = \{0, 1\}$ we only note---as already mentioned in example 3---that the (univariate) logistic regression model is equivalent to the association model \eqref{Verlobung}, so that no new aspects arise by using the latter model.

\section{Résumé and Discussion}
For a pair of random vectors $(X,Y)$ we have looked at semi-parametric association models with log-bilinear association---which include multivariate linear logistic regression, log-linear models for contingency tables as well as univariate and
multivariate linear regression models. Given a sample ($x_i, y_i),\ i= 1,\;\ldots,\;n$, the statistical inference for the odds-ratio parameter $\t$ (i.e. test and confidence regions) depends on the distribution of $\hat{\t}$ which typically is asymptotic normal and its covariance has to be estimated. The asymptotic approaches depend on the sampling scheme (conditional on $X$ resp. $Y$ or unconditional) and differ if $Y$ resp. $X$---or both---have finite support. We have shown however, that the estimated asymptotic covariance matrix of $\hat{\t}$ is invariant against the usual sampling schemes and does not depend on the support of $X$ or $Y$ being finite or arbitrary.
\par\bigskip\noindent
More precisely, we first considered the case where $X$ and $Y$ both have finite support. Then the log-bilinear odds-ratio model is a log-linear model for the expectations of the corresponding contingency table and by theorem \ref{Th3} the estimate $\hat{\Sigma}_{\hat{\t}}$ of the asymptotic covariance matrix $\Sigma_{\hat{\t}}$ is invariant against the common sampling schemes. Explicit representations for computing the matrix $\Sigma_{\hat{\t}}$ are given in theorem \ref{Th4}, \ref{Th5} and corollary \ref{Th6}. Allowing arbitrary support for $X$ but finite support for $Y$, the log-bilinear association model is a multivariate
linear logistic regression model. Our theorem \ref{Th7} implies that in this case the asymptotic covariance matrix of $\hat{\t}$ coincides with $\Sigma_{\hat{\t}}$ (where $X$ had finite support too).
\par\bigskip\noindent
To cover the general case with arbitrary supports of $X$ and $Y$ we looked at sampling conditional on $Y$ and an asymptotic approach where the set of conditioning values remains fixed. Combining the findings here with our earlier work we found that for a given sample the estimated asymptotic covariance matrix of $\hat{\t}$ coincides with the one computed for the observed contingency table under fixed cells asymptotics. And a dual result holds for sampling conditional on $X$ instead of $Y$.
\par\bigskip\noindent
Hence for asymptotic inference about the association parameter $\t$ one may assume any of the above sampling schemes and the statistical analysis of the sample can proceed as if both $X$ and $Y$ have finite support. Probably the most simple approach is to analyze the observed contingency table containing the counts $r_{jk}$ for all observed combinations of $x$-values and $y$-values using a log-linear model. Then an estimate of $\Sigma_{\hat{\t}}$ is obtained from corollary \ref{Th6} by using the estimate $\hat{D} = diag\{\vec{\hat{\mu}}\}$ instead of $D$. As a first application we have explained how our results allow to compute the asymptotic power for test of linear hypothesis about $\t$ and to determine the sample sizes to achieve a given power. Furthermore we have recovered the linear and log-linear regression model for univariate $Y$ from a more general log-bilinear association model.
\par\bigskip\noindent
Semiparametric association models do not restrict the marginal distributions of $X$ and $Y$. But more important, statistical inference about the association parameter $\t$ is possible for conditional sampling on either $X$ or $Y$. If $X$ is considered as an "input" and $Y$ as an "output" then sampling conditional on $X$ is a natural approach. However in certain situations sampling conditional on $Y$ may be advantageous, e.g. takes less time or money. For finite $Y$, for example, sampling conditional on $Y$ is very popular in epidemiology (case-control-studies) and econometrics (choice-based samples)---mainly because of their retrospective character. But as we have shown, sampling conditional on $Y$ may also be used if $Y$ has arbitrary support. In particular, using for univariate continuous $Y$ the more general log-bilinear odds-ratio model instead of the linear regression model, allows asymptotic inference even for the regression parameter when sampling is conditional on $Y$.
\par\bigskip\noindent
If sampling conditional on $X$ or $Y$ are an option, then the sampling scheme can be chosen to maximize the power of the test concerning the hypothesis of primary interest. This is well known for binary $X$ and $Y$ in the context of $2\times2$-tables and our results allow similar considerations for arbitrary $X$ and $Y$.

\section*{Appendix}
\addcontentsline{toc}{section}{Appendix}
In appendix \ref{A} and \ref{A.2} we summarize some definitions and results from linear algebra, which are used freely throughout the paper without explicit reference. In appendix \ref{C} the proofs of the theorems are given.
\begin{appendix}
\section{Inner Products and Orthogonal Projections}\label{A}
Any positive-definite symmetric $(I\times I)$ matrix $D$ induces an inner product on the vector space $\mathbb{R}^I$ given by $\left\langle a,b\right\rangle_D=a^TDb$, and orthogonality with respect to this inner product will be called $D$-orthogonality, denoted by $\bot_D$. 
\par\bigskip
\noindent
Consider a linear subspace $\mathscr{N}$ of $\mathbb{R}^I$ and a matrix $X$ whose columns form a basis of $\mathscr{N}$. The $D$-orthogonal projection $P^D_{\mathscr{N}}:\mathbb{R}^I\rightarrow \mathscr{N}$ onto $\mathscr{N}$ can be represented as an $I\times I$ matrix
\begin{equation}\label{Pr1}
P^D_{\mathscr{N}}=\ X(X^TDX)^{-1}X^TD.
\end{equation}
Some basic properties are
\begin{align}
&P^D_{\mathscr{N}}P^D_{\mathscr{N}}=P^D_{\mathscr{N}}\label{Pr2}\\
&\left(P^D_{\mathscr{N}}\right)^T=DP^D_{\mathscr{N}}D^{-1}\label{Pr3}\\
&\left(P^D_{\mathscr{N}}\right)^TDP^D_{\mathscr{N}}=DP^D_{\mathscr{N}}.\label{Pr10}
\end{align}
The $D$-orthogonal projection onto the $D$-orthogonal complement $\mathscr{N}^{\bot_D}= D^{-1}[N^{\bot}]$ satisfies
\begin{align}
P^D_{\mathscr{N}^{\bot_D}}&=D^{-1}P^{D^{-1}}_{\mathscr{N}^{\bot}}D\\
\mathbb{I}&=P^D_{\mathscr{N}}+P^D_{\mathscr{N}^{\bot_D}}\label{Pr4}
\end{align}
with the identity matrix $\mathbb{I}$. For another linear subspace $\mathscr{M}\subset\mathbb{R}^I$ it holds
\begin{align}
\mathscr{N}\subset\mathscr{M}\ &\Rightarrow\ P^D_{\mathscr{N}}P^D_{\mathscr{M}}=P^D_{\mathscr{N}}=P^D_{\mathscr{M}}P^D_{\mathscr{N}}\label{Pr7}\\
\mathscr{N}\bot_D\mathscr{M}\ &\Rightarrow\ P^D_{\mathscr{N}\oplus\mathscr{M}}=P^D_{\mathscr{N}}+P^D_{\mathscr{M}}.
\end{align}

\section{Kronecker Products}\label{A.2}
The Kronecker product of the two matrices, denoted by $A\otimes B$ is defined as the partitioned matrix (cf.\ Graham, 1981 \cite{graham})
\begin{equation}
A\otimes B= \begin{pmatrix} a_{11}B & a_{12}B & \ldots & a_{1n}B\\
													a_{21}B & a_{22}B & \ldots & a_{2n}B\\
													\vdots & \vdots &  & \vdots\\
													a_{m1}B & a_{m2}B & \ldots & a_{mn}B\end{pmatrix}.
\end{equation}
Some basic properties are
\begin{align}
& (\alpha A)\otimes(\beta B) = (\alpha\beta)(A\otimes B)\\
&(A +B) \otimes C = A \otimes C +B \otimes C\\
&A \otimes (B + C) = A \otimes B +A \otimes C\\
&(A \otimes B)^T = A^T \otimes B^T\\
&(A\otimes B)(C\otimes D) = AC\otimes BD\\
&(A \otimes B)^{-1}= A^{-1} \otimes B^{-1}\\
& (AYB)^{\vec{}}=(B^T \otimes A)\vec{Y}.
\end{align}

\section{Proofs}\label{C}

\subsection{Proof of Theorem \ref{Th3}}\label{Proof_Th3}
\eqref{100} restricts $\vec{\psi}\text{°}$ to the linear subspace $\mathscr{Q}=\{Z\text{°}\vec{\t}|\vec{\t}\in \mathbb{R}^L\}$ i.e.
\begin{equation}
\vec{\psi}\text{°} = Z\text{°}\vec{\t}
\end{equation}
where $Z\text{°}$ is assumed to have rank $L$. The parameters $\a$, $\r_j$, $\g_k$ and $\psi_{jk}$ are linear functions of the log-expectation $\eta$ and in particular
\begin{align}
\g_k&=\eta_{0k}-\eta_{00}\\
\psi_{jk}&=\eta_{jk}+\eta_{00} -\eta_{j0}-\eta_{0k}\label{psi}.
\end{align}
Then $\g\text{°}$ and $\psi\text{°}$ are given by
\begin{align}
\begin{pmatrix} \g\text{°} \\ \vec{\psi}\text{°} \end{pmatrix}&=\begin{pmatrix} B^T \\ C^T \end{pmatrix}\vec{\eta}. 
\end{align}
The columns of $B$ are orthogonal to the row space $\mathscr{R}$ and hence
\begin{equation}\label{B}
B^TP_{\mathscr{R}}^D=0.
\end{equation}
The columns of $C$\label{Q63} span the orthogonal complement of the marginal space $\mathscr{T}$ and thus
\begin{equation}\label{OR14}
P_{\mathscr{T}}^D D^{-1}C =0, \hspace{1cm} C^TP_{\mathscr{N}}^D = 0
\end{equation}
since $\mathscr{N}\subset\mathscr{T}$. The parameter $\l=(\g\text{°},\t)$ is linked to $(\g\text{°},\psi\text{°})$ in the following way
\begin{equation}
\begin{pmatrix} \g\text{°} \\ \vec{\psi}\text{°} \end{pmatrix}=\begin{pmatrix} \g\text{°} \\ Z\text{°}\vec{\t} \end{pmatrix}=
\begin{pmatrix} \mathbb{I}_K & 0 \\ 0 & Z\text{°} \end{pmatrix} \vec{\l}.
\end{equation}
Since the $(JK)\times L$\label{Q62} matrix $Z\text{°}$ has rank $L$ and a left inverse $Z\text{°}^{-}=(Z\text{°}^{T}Z\text{°})^{-1}Z\text{°}^{T}$ we get
\begin{equation}\label{C20}
\l=\begin{pmatrix} \mathbb{I}_K & 0 \\ 0 & Z\text{°}^{-} \end{pmatrix} \begin{pmatrix} \g\text{°} \\ \vec{\psi}\text{°} \end{pmatrix}
=\begin{pmatrix} \mathbb{I}_K & 0 \\ 0 & Z\text{°}^{-} \end{pmatrix} \begin{pmatrix} B^T \\ C^T \end{pmatrix}\vec{\eta}
=\begin{pmatrix} B^T\\Z\text{°}^{-}C^T \end{pmatrix}\vec{\eta}.
\end{equation}
Hence the asymptotic covariance matrix of the estimator $\hat{\l}$ can be derived from $\Sigma_{\hat{\eta}}$ as
\begin{equation}\label{Q55}
\Sigma_{\hat{\l}}=\begin{pmatrix} B^T\\Z\text{°}^{-}C^T \end{pmatrix} \Sigma_{\hat{\eta}} \begin{pmatrix} B,&CZ\text{°}^{-T} \end{pmatrix}.
\end{equation}
Using the block notation
\begin{equation}\label{AAA1}
\Sigma_{\hat{\l}}=\begin{pmatrix} \Sigma_{\hat{\g}\text{°}} & \Sigma_{\hat{\g}\text{°}\hat{\t}}\\
														 								 \Sigma_{\hat{\t}\hat{\g}\text{°}} & \Sigma_{\hat{\t}}	\end{pmatrix}
\end{equation}
we get the asymptotic covariance matrix of $\hat{\t}$ as
\begin{equation} \label{Sigma3}
\begin{split}
\Sigma_{\hat{\t}}&=Z\text{°}^{-}[C^T P_{\mathscr{H}}^D-C^TP_{\mathscr{N}}^D]D^{-1}CZ\text{°}^{-T}\\
																		&=Z\text{°}^{-}C^T P_{\mathscr{H}}^D D^{-1}CZ\text{°}^{-T}.
\end{split}
\end{equation}
\begin{flushright}
$\square$
\end{flushright}

\subsection{Proof of Theorem \ref{Th4}}\label{Proof_Th4}
The space of the log-expectation $\mathscr{H}$ can then be decomposed into the direct sum
\begin{equation}\label{AA2}
\begin{split}
\mathscr{H}=\mathscr{T} \oplus \mathscr{Z}'\hspace{1cm}\text{with}\hspace{1cm}\mathscr{Z}'=\mathscr{H}\cap\mathscr{T}^{\bot_D}.
\end{split}
\end{equation}
$\mathscr{Z}'$ is the $D$-orthogonal complement of $\mathscr{T}$ in $\mathscr{H}$. Applying the orthogonal projection $P_{\mathscr{T}^{\bot_D}}^D$ on $\mathscr{H}$ yields 
\begin{equation}
P_{\mathscr{T}^{\bot_D}}^D [\mathscr{H}]
=\mathscr{Z}'
\end{equation}
and hence the columns of the $(I\times L)$-matrix 
\begin{equation}\label{OR1}
V=P_{\mathscr{T}^{\bot_D}}^D Z
\end{equation}
span $\mathscr{Z}'$. Using the representation
\begin{equation}\label{1001}
P_{\mathscr{Z}'}^D=V(V^TDV)^{-1}V^TD
\end{equation}
and $\eqref{OR14}$ we get
\begin{equation}\label{OR11}
C^T P_{\mathscr{H}}^D D^{-1}C\ =\ C^T P_{\mathscr{T}}^D D^{-1}C + C^T P_{\mathscr{Z}'}^D D^{-1}C\ =\ C^T V(V^TDV)^{-1}V^T C.
\end{equation}
Since the columns of $C$ are elements of $\mathscr{T}^{\bot}$ we get
\begin{equation}
\begin{split}
Z^TC=\;Z^TP_{\mathscr{T}^{\bot}}^{D^{-1}}C=\;Z^T(P_{\mathscr{T}^{\bot_D}}^D)^T C=V^TC
\end{split}
\end{equation}
For $j,k>0$ the $jk$-th row of $C^TZ$
\begin{equation}
(C^TZ)_{jk}=c_{jk}^TZ=\vec{e}_{jk}^{\;T}Z + \vec{e}_{00}^{\;T}Z - \vec{e}_{j0}^{\;T}Z - \vec{e}_{0k}^{\;T}Z=z_{jk}
\end{equation}
holds because $z_{j0}=z_{0k}=0$ and therefore
\begin{equation}\label{OR12}
C^TV=C^TZ=Z\text{°}.
\end{equation}
With $\eqref{OR11}$ and \eqref{OR1} this leads to different representations of $\Sigma_{\hat{\t}}$ from $\eqref{Sigma3}$ 
\begin{equation}\label{OR22}
\begin{split}
\Sigma_{\hat{\t}}=\ \ &Z\text{°}^-(C^T V(V^TDV)^{-1}V^T C)Z\text{°}^{-T}\\
													 =\ \ &(V^TDV)^{-1}\\
													 =\ \ &(Z^T(P_{\mathscr{T}^{\bot_D}}^D)^TDP_{\mathscr{T}^{\bot_D}}^DZ)^{-1}\\
													 =\ \ &(Z^TP_{\mathscr{T}^{\bot}}^{D^{-1}}DZ)^{-1}\\
													 =\ \ &(Z^TC(C^TD^{-1}C)^{-1}C^TD^{-1}DZ)^{-1}
													 \end{split}
\end{equation}
and thus with $\eqref{OR12}$ the representation \eqref{OR2} is obtained.$\hspace{5.8cm}\square$

\subsection{Proof of Theorem \ref{Th5}}\label{Proof_Th5}
The marginal space $\mathscr{T}$ can be decomposed into $\mathscr{T}=\ \mathscr{R}\oplus \mathscr{C}'$ where $\mathscr{C}'$ is spanned by the columns of the matrix $E=\ (\vec{e}_{+1},\;\ldots,\;\vec{e}_{+K})$. The columns of the $(I\times K)$ matrix
\begin{equation}\label{Q90}
S=P_{\mathscr{R}^{\bot_D}}^D E
\end{equation}
span the $D$-orthogonal complement $\mathscr{C}''=\mathscr{T}\cap\mathscr{R}^{\bot_D}$ of $\mathscr{R}$ within $\mathscr{T}$. The direct decomposition 
\begin{equation}\label{AA3}
\mathscr{T}=\mathscr{R}\oplus \mathscr{C}''
\end{equation}
leads to
\begin{equation}\label{OR5}
Z^TDP_{\mathscr{T}^{\bot_D}}^DZ=\ \;Z^TDZ-Z^TDP_{\mathscr{T}}^DZ
															=Z^TDP_{\mathscr{R}^{\bot_D}}^DZ-Z^TDP_{\mathscr{C}''}^DZ.
\end{equation}
Since $\Sigma_{\hat{\t}}$ is invariant against the underlying distribution model we will assume for the rest of the proof the product multinomial sampling for rows.  The rows are independent of each other and we know that $Cov_{\text{\tiny{MR}}}(\vec{R})=DP_{\mathscr{R}^{\bot_D}}^D$ (cf.\ \eqref{(3.14)}) where the index MR refers to the sampling scheme.\label{OR6} The $D$-orthogonal projection on $\mathscr{C}''$can now be specified as
\begin{equation}\label{OR7}
\begin{split}
DP_{\mathscr{C}''}^D=\ \;&DS(S^TDS)^{-1}S^TD\\
										=\ \;&DP_{\mathscr{R}^{\bot_D}}^D E(E^TDP_{\mathscr{R}^{\bot_D}}^D E)^{-1}(DP_{\mathscr{R}^{\bot_D}}^D E)^T\\
										=\ \;&Cov_{\text{\tiny{MR}}}(\vec{R})E(E^TCov_{\text{\tiny{MR}}}(\vec{R})E)^{-1}E^TCov_{\text{\tiny{MR}}}(\vec{R}).
\end{split}
\end{equation}
Together with $\eqref{OR5}$ we obtain from $\eqref{OR2}$ 
\begin{equation}\label{OR9}
\begin{split}
\Sigma_{\hat{\t}}\ =&\left[Z^TDP_{\mathscr{T}^{\bot_D}}^DZ\right]^{-1}\\
											 =&\left[Z^TDP_{\mathscr{R}^{\bot_D}}^DZ-Z^TDP_{\mathscr{C}''}^DZ\right]^{-1}\\
											 =&\left[Z^T Cov_{\text{\tiny{MR}}}(\vec{R})Z-Z^TCov_{\text{\tiny{MR}}}(\vec{R})E(E^TCov_{\text{\tiny{MR}}}(\vec{R})E)^{-1}E^TCov_{\text{\tiny{MR}}}(\vec{R})Z\right]^{-1}.
\end{split}
\end{equation}
\begin{flushright}
$\square$
\end{flushright}

\subsection{Proof of Theorem \ref{Th7}}\label{Proof_Th7}
We first determine the matrices $Cov(U(\l))$ and $\Sigma_{\hat{\l}}$ with $E=(\vec{e}_{+1},\;\ldots,\;\vec{e}_{+K})$ and $Z = (z_{jk})_{jk}$, $z_{jk}$ being the rows of $Z$. The score vector may be written as
\begin{equation}\label{Q29}
U(\l)=\begin{pmatrix} U_{\g}(\l) \\ U_{\t}(\l)\end{pmatrix}
					=\begin{pmatrix} \left(\sum_{j=0}^J (R_{jk}-n_j p^X_{jk})\right)_{k=1,\;\ldots,\;K} \\ 
														\sum_{j=0}^J \sum_{k=1}^K (R_{jk}-n_j p^X_{jk})z_{jk}^T\end{pmatrix}
					=\begin{pmatrix} E^T(\vec{R}-\vec{\mu}) \\ 
														Z^T(\vec{R}-\vec{\mu})\end{pmatrix}.
\end{equation}
Hence
\begin{equation}\label{I2}
Cov(U(\l))=
\begin{pmatrix} E^T Cov(\vec{R})E & E^T Cov(\vec{R})Z \\ Z^T Cov(\vec{R})E & Z^T Cov(\vec{R})Z \end{pmatrix}.
\end{equation}
The matrix $\Sigma_{\hat{\l}}$ has a block representation \eqref{AA1} and we know from \eqref{Q5}, \eqref{B} and \eqref{Sigma3} that
\begin{align}
\Sigma_{\hat{\l}}=\ \;&\begin{pmatrix} \Sigma_{\hat{\g}\text{°}} & \Sigma_{\hat{\g}\text{°}\hat{\t}}\\
														 								 \Sigma_{\hat{\t}\hat{\g}\text{°}} & \Sigma_{\hat{\t}}	\end{pmatrix}\notag\\
														=\ \;&\begin{pmatrix} [B^T P_{\mathscr{H}}^D-B^TP_{\mathscr{R}}^D]D^{-1}B
																						& [B^T P_{\mathscr{H}}^D-B^TP_{\mathscr{R}}^D]D^{-1}CZ\text{°}^{-T}\\
														 								 ([B^T P_{\mathscr{H}}^D-B^TP_{\mathscr{R}}^D]D^{-1}CZ\text{°}^{-T})^T  																																		     &	Z\text{°}^-C^T P_{\mathscr{H}}^D D^{-1}CZ\text{°}^{-T}  \end{pmatrix}\\
														 								 =\ \;&\begin{pmatrix} B^T P_{\mathscr{H}}^DD^{-1}B
																						& B^T P_{\mathscr{H}}^DD^{-1}CZ\text{°}^{-T}\\
														 								 (B^T P_{\mathscr{H}}^DD^{-1}CZ\text{°}^{-T})^T  																																		     &	Z\text{°}^-C^T P_{\mathscr{H}}^D D^{-1}CZ\text{°}^{-T}  \end{pmatrix}.\notag
\end{align}
Each of these blocks will be determined separately similarly to \eqref{OR11}. First we obtain $B^TE=\mathbb{I}_K$ since
\begin{align*}
(B^TE)_{kk}&=b_k^T\vec{e}_{+k}=\vec{e}_{0k}^{\;T}\vec{e}_{+k} - \vec{e}_{00}^{\;T}\vec{e}_{+k}=\vec{e}_{0k}^{\;T}\vec{e}_{+k}=1 &\text{, for $k>0$.}\\
(B^TE)_{ki}&=b_k^T\vec{e}_{+i}=\vec{e}_{0k}^{\;T}\vec{e}_{+i} - \vec{e}_{00}^{\;T}\vec{e}_{+i}=\vec{e}_{0k}^{\;T}\vec{e}_{+i}=0 &\text{, for $i>0$ and $i\neq k$.}
\end{align*}
It follows from \eqref{B}
\begin{align}\label{Q7}
B^T D^{-1}Cov(\vec{R})E=B^T(\mathbb{I}- P_{\mathscr{R}}^D)E 																																										=\;B^TE-B^TP_{\mathscr{R}}^D E=\mathbb{I}_K.
\end{align}
The D-orthogonal decomposition \eqref{AA2} as well as \eqref{OR14}, \eqref{1001} and \eqref{OR12} lead to 
\begin{align*} 
\Sigma_{\hat{\g}\text{°},\hat{\t}}\ =\ \;&B^T P_{\mathscr{T}}^D D^{-1}C Z\text{°}^{-T}+B^T P_{\mathscr{Z}'}^D D^{-1}CZ\text{°}^{-T}=\ \;B^T P_{\mathscr{Z}'}^D D^{-1}CZ\text{°}^{-T}\\
				   						=\ \;&B^TV(V^TDV)^{-1}V^TDD^{-1}C Z\text{°}^{-T}
				   						=\ \;B^TV(V^TDV)^{-1}.
\end{align*}
And the $D$-orthogonal decomposition $\mathscr{T}=\mathscr{R}\oplus \mathscr{C}''$ \eqref{AA3} together with \eqref{B}, \eqref{1001}, \eqref{OR7} and \eqref{Q7} yields
\begin{align*}
\Sigma_{\hat{\g}\text{°}}\ =\ \;&B^T P_{\mathscr{T}}^D D^{-1}B + B^T P_{\mathscr{Z}'}^D D^{-1}B
    =\ \;B^T P_{\mathscr{R}}^D D^{-1}B+ B^T P_{\mathscr{C}''}^D D^{-1}B + B^T P_{\mathscr{Z}'}^D D^{-1}B\\
		=\ \;&B^T P_{\mathscr{C}''}^D D^{-1}B + B^T P_{\mathscr{Z}'}^D D^{-1}B
		=\ \;B^T P_{\mathscr{C}''}^D D^{-1}B + B^T V(V^TDV)^{-1}V^TB\\
		=\ \;&B^T D^{-1} (Cov(\vec{R})E(E^TCov(\vec{R})E)^{-1}E^TCov(\vec{R}))D^{-1}B + B^T V(V^TDV)^{-1}V^TB\\
		=\ \;&(E^TCov(\vec{R})E)^{-1}+ B^T V(V^TDV)^{-1}V^TB.
\end{align*}
Using \eqref{OR2} we can summarize this into
\begin{equation}\label{Q8}
\Sigma_{\hat{\l}}=\begin{pmatrix} (E^TCov(\vec{R})E)^{-1}+ B^T V(V^TDV)^{-1}V^TB & B^TV(V^TDV)^{-1} \\
																							(V^TDV)^{-1}V^TB 												& (V^TDV)^{-1} \end{pmatrix}.
\end{equation}
To prove the theorem we further examine \eqref{Q8}. The term $(V^TDV)^{-1}$ is known from previous considerations. We now have a closer look at the remaining term $V^TB$. Since the first $K+1$ rows of $Z$ are equal to zero and the first $K+1$ rows of $B$ are the only rows of $B$  with entries non-equal to zero we get $Z^TB=\ 0$. From \eqref{OR1}, \eqref{B}, \eqref{OR7} and \eqref{Q7} it follows
\begin{equation}
\begin{split}
V^TB=\ \; &Z^T DP_{\mathscr{T}^{\bot_D}}^D D^{-1}B
													 =\ \;Z^TDD^{-1}B-Z^TDP_{\mathscr{T}}^DD^{-1}B\\
	 												 =\ \;&-Z^TD\left(P_{\mathscr{R}}^D + P_{\mathscr{C}''}^D\right)D^{-1}B
	  											 =\ \;-Z^T(B^TP_{\mathscr{R}}^D)^T - Z^TDP_{\mathscr{C}''}^DD^{-1}B\label{Q9}\\
	  											 =\ \;&-Z^TDP_{\mathscr{C}''}^DD^{-1}B
	  											 =\ \;-Z^TCov(\vec{R})E(E^TCov(\vec{R})E)^{-1}E^TCov(\vec{R}))D^{-1}B\\
	  											 =\ \;&-Z^TCov(\vec{R})E(E^TCov(\vec{R})E)^{-1}.
\end{split}
\end{equation}
After determining all components of $\Sigma_{\hat{\l}}$ we are going to invert $Cov(U(\l))$ using \eqref{Matrix1}. 
For
\begin{equation*}
A=Cov(U(\l))=\begin{pmatrix} E^T Cov(\vec{R})E & E^T Cov(\vec{R})Z \\ Z^T Cov(\vec{R})E & Z^T Cov(\vec{R})Z \end{pmatrix}
\end{equation*}
we compute $A^{-1}=Cov(U(\l))^{-1}$ and let
\begin{equation*}
L=\;E^T Cov(\vec{R})E,\hspace{1cm}
M=\;E^T Cov(\vec{R})Z,\hspace{1cm}
G=\;M^T ,\hspace{1cm}
H=\;Z^T Cov(\vec{R})Z.
\end{equation*}
Then \eqref{OR9}, \eqref{OR2} and \eqref{Q9} lead to
\begin{equation*}
\begin{split}
N=\ \ &H-GL^{-1}M
=\;Z^T Cov(\vec{R})Z-Z^T Cov(\vec{R})E(E^T Cov(\vec{R})E)^{-1}E^T Cov(\vec{R})Z\\
=\ \ &\Sigma_{\hat{\t}}^{-1}
=\;V^TDV\\
-L^{-1}M=\ \ &-(E^T Cov(\vec{R})E)^{-1}E^T Cov(\vec{R})Z
=\;B^TV\\
-GL^{-1}=\ \ &-Z^T Cov(\vec{R})E(E^T Cov(\vec{R})E)^{-1}
=\;V^TB
\end{split}
\end{equation*}
and accordingly
\begin{equation*}
L^{-1}+L^{-1}MN^{-1}GL^{-1}=\ (E^T Cov(\vec{R})E)^{-1}+B^TV(V^TDV)^{-1}V^TB.
\end{equation*}
Summing up, \eqref{Matrix1} and \eqref{Q8} yields
\begin{equation*}
\begin{split}
Cov(U(\l))^{-1}=\ \;&\begin{pmatrix} (E^TCov(\vec{R})E)^{-1}+ B^T V(V^TDV)^{-1}V^TB & B^TV(V^TDV)^{-1} \\(V^TDV)^{-1}V^TB & (V^TDV)^{-1} \end{pmatrix}
										=\;\Sigma_{\hat{\l}}.
\end{split}
\end{equation*}
\begin{flushright}
$\square$
\end{flushright}
\end{appendix}

\thispagestyle{plain}

\vspace{1cm}
\noindent
\textit{Date: 11 April 2012\\
This paper is available for download at:} http://www.math.uni-bremen.de/$\sim$osius
\end{document}